\documentclass[10pt,leqno]{article}
\usepackage{amsfonts,amssymb,amsmath,amsthm,latexsym}
\usepackage{color,graphicx}

\newcommand{\CC}{\mathbb{C}}

\newtheorem{teor}{Theorem}
\newtheorem{lema}{Lemma }
\newtheorem{prop}{Proposition}
\newtheorem{obs}{Remark }

\newcommand{\ra}{\rightarrow}

\newcommand{\ctd}{$\blacksquare$}
\newcommand{\demo}{\noindent {\it Proof:   }}

\newcommand{\xe}{x_{\varepsilon}}
\newcommand{\gv}{\gamma_{\varepsilon}}
\newcommand{\al} {\alpha}
\newcommand{\ve} {\varepsilon}

\newcommand{\ds} {\displaystyle}

\newcommand{\lt} {L^{2}(-\frac{\widetilde{T}}{2},\frac{\widetilde{T}}{2})}

\newcommand{\lm} {\lambda}

\newcommand{\vf} {\varphi_{\ve}^{-1}}

\numberwithin{equation}{section}
\begin{document}
\title{A singular controllability problem with vanishing viscosity}

\author{Ioan Florin Bugariu \thanks{Facultatea de Stiinte Exacte, Universitatea din
Craiova, 200585, Romania. E-mail: {\tt
florin$\_$bugariu$\_$86@yahoo.com}} \and Sorin Micu
\thanks{Facultatea de Stiinte Exacte, Universitatea din
Craiova, 200585, Romania. E-mail: {\tt sd$\_$micu@yahoo.com}.}}

\maketitle

\begin{abstract}
The aim of this paper is to answer the question: Do the controls
of a vanishing viscosity approximation of the one dimensional
linear wave equation converge to a control of the conservative
limit equation? The characteristic of our viscous term is that it
contains the fractional power $\alpha$ of the Dirichlet Laplace
operator. Through the parameter $\alpha$ we may increase or
decrease the strength of the high frequencies damping which allows
us to cover a large class of dissipative mechanisms. The viscous
term, being multiplied by a small parameter $\ve$ devoted to tend
to zero, vanishes in the limit. Our analysis, based on moment
problems and biorthogonal sequences, enables us to evaluate the
magnitude of the controls needed for each eigenmode and to show
their uniform boundedness with respect to $\ve$, under the
assumption that $\alpha\in [0,1)\setminus
\left\{\frac{1}{2}\right\}$. It follows that, under this
assumption, our starting question has a positive answer.

\


\

\noindent{{\it Mathematical subject codes:} 93B05, 30E05, 35E20.}
\end{abstract}

\section{Introduction}

For $T>0$, we consider the one-dimensional linear wave equation with
``lumped" control
\begin{equation}\label{ec.in0}
\left\{
\begin{array}{ll}
w_{tt}(t,x)-\partial_{xx}^2 w(t,x)=v(t)f(x)&(t,x)\in(0,T)\times(0,\pi)\\
w(t,0)=w(t,\pi)=0& t\in(0,T)\\
w(0,x)=w^0(x),\,\, w_t(0,x)=w^1(x)&x\in(0,\pi),
\end{array}
\right.
\end{equation}
where the profile $f\in L^2(0,\pi)$ is given and verifies
$\ds\widehat{f}_n\neq0$ for every $n\geq 1$. Here and in the
sequel, given any function $g\in L^2(0,\pi)$, we denote by
$\widehat{g}_n$ the $n-$th Fourier coefficient of $g$,
$$\widehat{g}_n=\int_{0}^{\pi}g(x)\sin(nx)dx\qquad (n\geq 1).$$

Equation \eqref{ec.in0} is said to be {\it null-controllable in time
$T>0$} if, for every initial data $(w^0,w^1)\in {\cal H}_0\subset
H^1_0(0,\pi)\times L^2(0,\pi)$, there exists a control $v\in
L^2(0,T)$ such that the corresponding solution of \eqref{ec.in0}
verifies
\begin{equation}\label{ec.fin}
w(T,\,\cdot\,)=w_t(T,\,\cdot\,)=0,
\end{equation}
where the space ${\cal H}_0$ is defined as follows
\begin{equation}\label{sh}{\cal H}_0=\left\{(w^0,w^1)\in
H_0^1(0,\pi)\times L^2(0,\pi)\ \left| \ \sum_{n\geq 1}
\frac{n^2\left|\widehat{w}_n^0\right|^2 +
\left|\widehat{w}_n^1\right|^2}{\left|\widehat{f}_n\right|^2}<\infty
\right\}\right. .\end{equation}

The controllability properties of \eqref{ec.in0} are by now
well-known (see, for instance, the monographs \cite{C,Tuc}). One
of the oldest methods used to study such controllability problems
consists in reducing them to a {\it moment problem} whose solution
is given in terms of an explicit biorthogonal sequence to a family
$\Lambda$ of exponential functions. For instance, this method was
used by Fattorini and Russell in the pioneering articles
\cite{Fat1,Fat2} to prove the controllability of the one
dimensional heat equation. In their case, the family $\Lambda$ has
only real exponential functions. On the contrary, when equation
\eqref{ec.in0} comes into discussion, the family $\Lambda$ is
given by $\left(e^{\mu_n t}\right)_{n\in\mathbb{Z}^*}$, where
$\mu_n=in$, $n\in\mathbb{Z}^*$, are the eigenvalues of the wave
operator $\left(\begin{array}{cc}0 & -I
\\-\partial_{xx}^2 & 0\end{array}\right)$ and are purely imaginary.
It follows easily that, \eqref{ec.in0} is null-controllable in
time $T$ if, and only if, for every initial data $(w^0,w^1)\in
{\cal H}_0$, there exists a solution $v\in L^2(0,T)$ of the
following moment problem:
\begin{equation}\label{pm0}
\int_{-\frac{T}{2}}^{\frac{T}{2}}v\left(t+\frac{T}{2}\right)e^{\overline{\mu}_nt}dt=
-\frac{e^{-\frac{T}{2}\overline{\mu}_n}}{\widehat{f}_{|n|}}
\left(\widehat{w}^1_{|n|}+\mu_n
\widehat{w}^0_{|n|}\right)\qquad(n\in \mathbb{Z}^*).
\end{equation}

In order to fix some ideas and to illustrate the method used in
this paper, let us briefly show how do we obtain a solution of
\eqref{pm0}. We begin by defining the function
\begin{equation}\label{defpsi}
\widetilde{\Psi}_m(z)=\frac{\sin(\pi(z+m))}{\pi(z+m)},
\end{equation}
which is an entire function of exponential type $\pi$ such that
$\ds\int_{\mathbb{R}}\left|\widetilde{\Psi}_m(x)\right|^2dx<\infty.$
It results from $Paley-Wiener$ Theorem that the Fourier transform
of $\widetilde{\Psi}_m$,
\begin{equation}
\widetilde{\theta}_m(t)=\frac{1}{2\pi}\int_{\mathbb{R}}\widetilde{\Psi}_m(x)e^{-ixt}dx\qquad
(m\in\mathbb{Z}^*),
\end{equation}
belongs to $L^2(-\pi,\pi)$. Moreover, from the inversion formula,
it follows that $\left(\widetilde{\theta}_m\right)_{m\in
\mathbb{Z}^*}$ forms a biorthogonal sequence to the family
$\Lambda=\left(e^{\mu_n t}\right)_{n\in \mathbb{Z}^*}$, i. e.
verify
\begin{equation}\label{bio}
\int_{-\pi}^{\pi} \widetilde{\theta}_m(t) e^{\overline{\mu}_n t}\,
dt = \delta_{mn}\qquad (m,n\in\mathbb{Z}^*).
\end{equation}

From the above properties, we deduce that a formal solution of the
moment problem \eqref{pm0} is given by
\begin{equation}\label{fsol0}
v(t)=-\sum_{m\in\mathbb{Z}^*}
\frac{e^{-\pi\overline{\mu}_m}}{\widehat{f}_{|m|}}
\left(\widehat{w}^1_{|m|}+\mu_m \widehat{w}^0_{|m|}\right)
\widetilde{\theta}_m\left(t-\pi \right)\qquad(t\in(0,2\pi)).
\end{equation}
In fact \eqref{fsol0} gives a true solution of \eqref{pm0} if the
right hand side of \eqref{fsol0} converges in $L^2(0,2\pi)$. For
each $(w^0,w^1)\in{\cal H}_0$, the convergence of this series
follows from the existence of a constant $C>0$ such that
\begin{equation}\label{bot}
\|\widetilde{\theta}_m\|_{L^2(-\pi,\pi)}\leq
C\qquad(m\in\mathbb{Z}^*),
\end{equation}
which is a consequence of the uniform boundedness (in $m$) of the
$L^2(\mathbb{R})-$norms of
$(\widetilde{\Psi}_m)_{m\in\mathbb{Z}^*}$ and Plancherel's
Theorem. Hence, for any initial data $(w^0,w^1)\in {\cal H}_0$,
the moment problem has at least a solution $v\in L^2(0,2\pi)$,
given by \eqref{fsol0}, and the controllability of \eqref{ec.in0}
in time $T=2\pi$ follows.

In many applications it is of interest to study the uniform
controllability properties of \eqref{ec.in0} when a viscous term
is introduced in the equation. Indeed, the mechanism of vanishing
viscosity is a common tool in the study of Cauchy problems or in
improving convergence of numerical schemes for hyperbolic
conservation laws and shocks. For instance, in \cite{IZ1,IZ2}, it
is proved that, by adding an extra numerical viscosity term, the
dispersive properties of the finite difference scheme for the
nonlinear Schr\"odinger equation become uniform when the mesh-size
tends to zero. This scheme reproduces at the discrete level the
properties of the continuous Schr\"odinger equation by dissipating
the high frequency numerical spurious solutions. On the other
hand, a viscosity term is introduced in \cite{Dip} to prove the
existence of solutions of hyperbolic equations. In both examples
the viscosity is devoted to tend to zero in order to recover the
original system. Thus, a legitimate question is related to the
behavior and the sensitivity of the controls during this process.
For instance, given $T>0$ and $\ve\in(0,1)$, one could consider
the perturbed wave equation
\begin{equation}\label{ec.in1}
\left\{
\begin{array}{ll}
u_{tt}(t,x)-\partial^2_{xx}u(t,x)+2\varepsilon
(-\partial_{xx}^2)^\alpha
u_t(t,x)=v_\varepsilon(t)f(x)&(t,x)\in(0,T)
\times(0,\pi)\\
u(t,0)=u(t,\pi)=0& t\in(0,T)\\
u(0,x)=u^0(x),\,\, u_t(0,x)=u^1(x)&x\in(0,\pi)
\end{array}
\right.
\end{equation}
and study the possibility of obtaining a control $v$ of
\eqref{ec.in0} as limit of controls $v_\ve\in L^2(0,T)$ of
\eqref{ec.in1}. Here and in what follows
$(-\partial^2_{xx})^\alpha $ denotes the fractional power of order
$\alpha\geq 0$ of the Dirichlet Laplace operator in $(0,\pi)$.
More precisely,
\begin{equation}
\label{eq.0oplap}
\begin{array}{c}

\vspace{2mm}

(-\partial^2_{xx})^\alpha :D((-\partial^2_{xx})^\alpha)\subset
L^2(0,\pi)\rightarrow
L^2(0,\pi),\\

\vspace{2mm}

\displaystyle \quad D((-\partial^2_{xx})^\alpha)=\left\{u\in
L^2(0,\pi)\,:\, u=\sum_{n\geq 1} a_n \sin(n x)\mbox{ and
}\sum_{n\geq 1} |a_n|^2 n^{4\alpha}<\infty
\right\},\\

\displaystyle u(x)=\sum_{n\geq 1} a_n \sin(n x) \,\longrightarrow
\, (-\partial^2_{xx})^\alpha u(x)=\sum_{n\geq 1} a_n n^{2\alpha}
\sin(n x).\end{array}\end{equation}

Equation \eqref{ec.in1} is dissipative and it can be easily
checked that, if $f=0$,
\begin{equation}\label{ener1}
\frac{d}{dt}\left(\|u(t)\|_{H_0^1}^2 + \|u(t)\|_{L^2}^2
\right)=-2\ve \int_0^\pi \left|(-\partial_{xx}^2)^\frac{\alpha}{2}
u_t(t,x)\right|^2 dt \leq 0.
\end{equation}
Hence, $2\varepsilon (-\partial_{xx}^2)^\alpha u_t(t,x)$
represents an added viscous term devoted to vanish as $\ve$ tends
to zero. However, the controllability properties of \eqref{ec.in1}
are poor. Indeed, the family of exponential functions
corresponding to this case is given  by $\Lambda=\left(e^{\nu_n
t}\right)_{n\in\mathbb{Z}^*}$, where $\nu_n=\ve
|n|^{2\alpha}+\mbox{sgn\,(n)}\sqrt{|n|^{4\alpha}-n^2}$. If
$\alpha>\frac{1}{2}$, we have that
$$\lim_{n\rightarrow -\infty}\nu_n=0,$$
which implies that the family $\Lambda$ is not minimal.
Consequently, equation \eqref{ec.in1} is not spectrally
controllable if $\alpha>\frac{1}{2}$.

Since we want to allow stronger dissipative terms which correspond
to the case $\alpha>\frac{1}{2}$, we perturb the wave equation
\eqref{ec.in0} in the following slightly different way
\begin{equation}\label{ec.in}
\left\{
\begin{array}{ll}
u_{tt}(t,x)-\partial^2_{xx}u(t,x)+2\varepsilon
(-\partial_{xx}^2)^\alpha u_t(t,x)+ \varepsilon^2
(-\partial_{xx}^2)^{2\alpha}
u(t,x)=v_\varepsilon(t)f(x)&(t,x)\in(0,T)
\times(0,\pi)\\
u(t,0)=u(t,\pi)=0& t\in(0,T)\\
u(0,x)=u^0(x),\,\, u_t(0,x)=u^1(x)&x\in(0,\pi).
\end{array}
\right.
\end{equation}
Equation \eqref{ec.in} is still dissipative. Indeed, if $f=0$, we
have that
\begin{equation}\label{ener}
\frac{d}{dt}\left(\|u(t)\|_{H_0^1}^2 + \ve^2 \|(-\partial
_{xx})^{\alpha}u(t)\|_{L^2}^2+ \|u(t)\|_{L^2}^2 \right)=-2\ve
\int_0^\pi \left|(-\partial_{xx}^2)^\frac{\alpha}{2}
u_t(t,x)\right|^2 dt \leq 0.
\end{equation}

Note that, if $\alpha\leq \frac{1}{2}$, the norm
$\sqrt{\|u\|_{H_0^1}^2 + \ve^2 \|(-\partial
_{xx})^{\alpha}u\|_{L^2}^2+ \|u\|_{L^2}^2 }$ is equivalent to
$\sqrt{\|u\|_{H_0^1}^2 + \|u\|_{L^2}^2}$ and the controllability
properties of \eqref{ec.in1} and \eqref{ec.in} are similar.
However, in the case $\alpha>\frac{1}{2}$, the controllability
properties of \eqref{ec.in} are better than those of
\eqref{ec.in1}. Hence, the term $\varepsilon^2
(-\partial_{xx})^{2\alpha} u(t,x)$ allows us to consider stronger
dissipation and also to simplify some of our estimates.

The aim of this paper is to study the controllability properties
of \eqref{ec.in} and their relation with the ones of
\eqref{ec.in0}. The controllability of \eqref{ec.in} is defined in
a similar way as for \eqref{ec.in0}. More precisely, given $T>0$
and $f\in L^2(0,\pi)$ with $\widehat{f}_n\neq0$ for $n\geq1$,
equation \eqref{ec.in} is {\em null-controllable in time $T$} if,
for any $(u^0,u^1)\in {\cal H}_0$, there exists a control
$v_\ve\in L^2(0,T)$ such that the corresponding solution $(u,u_t)$
of \eqref{ec.in} verifies
\begin{equation}\label{ec.coco}
u(T,\,\cdot\,)=u_t(T,\,\cdot\,)=0.
\end{equation}
The null-controllability problem is equivalent to find, for every
initial data $(u^0,u^1)\in {\cal H}_0$, a solution $v_\varepsilon
\in L^2(0,T)$ of the following moment problem:
\begin{equation}\label{pm}
\int_{-\frac{T}{2}}^{\frac{T}{2}}v_{\ve}\left(t+\frac{T}{2}\right)e^{\overline{\lm}_{n}t}dt=-
 \frac{e^{-\overline{\lm}_{n}\frac{T}{2}}}{\widehat{f}_{|n|}}
 \left(\widehat{u}^1_{|n|}+\lm_n\widehat{u}^0_{|n|}\right)
 \qquad(n\in\mathbb{Z}^*),
\end{equation}
where $\lm_n=in+\varepsilon |n|^{2\alpha},$ are the eigenvalues of
the operator $\left(\begin{array}{cc}0 & -I
\\-\partial_{xx}^2 +\varepsilon^2 (-\partial_{xx}^2)^{2\alpha}&
2\ve(-\partial_{xx}^2)^\alpha\end{array}\right)$ corresponding to the ``adjoint" problem of
\eqref{ec.in}.

As in \eqref{pm0}, if we have at our disposal a biorthogonal
sequence to the family $\left(e^{\lm_n
t}\right)_{n\in\mathbb{Z}^*}$, denoted by
$\left(\theta_m\right)_{m\in\mathbb{Z}^*}$, we can give
immediately a formal solution of \eqref{pm},
\begin{equation}\label{fsol}
v_\varepsilon(t)=-\sum_{m\in\mathbb{Z}^*}
 \frac{e^{-\overline{\lm}_{m}\frac{T}{2}}}{\widehat{f}_{|m|}}
 \left(\widehat{u}^1_{|m|}+\lm_m\widehat{u}^0_{|m|}\right)\theta_{m}\left(t-\frac{T}{2}\right) \qquad
(t\in(0,T)).
\end{equation}

This time the family $\left(e^{\lm_n t}\right)_{n\in
\mathbb{Z}^*}$ has no longer purely imaginarily exponents like in
\eqref{ec.in0}. Thus, it is not so easy as for \eqref{ec.in0} to
give explicit entire functions $(\Psi_m)_{m\in\mathbb{Z}^*}$ whose
Fourier transforms define a biorthogonal sequence
$(\theta_m)_{m\in\mathbb{Z}^*}$ to  $\left(e^{\lm_n
t}\right)_{n\in\mathbb{Z}^*}$. Moreover, we cannot guarantee
anymore the boundedness of the sequence
$(\theta_m)_{m\in\mathbb{Z}^*}$ and \eqref{bot} will be replaced
by an estimate of the form
\begin{equation}\label{ein} ||\theta_m||_{L^2}\leq C e^{\beta |\Re
(\lambda_m)|}\qquad (m\in\mathbb{Z}^*).\end{equation} Note that
$||\theta_m||_{L^2}$ may become exponentially large as $m$ goes to
infinity. By taking into account the damping mechanism introduced
in equation \eqref{ec.in}, this growth estimate guarantees the
convergence of series \eqref{fsol} for each initial data
$(u^0,u^1)\in{\cal H}_0$, if $T$ is large enough. However, in
order show that a control time $T$ independent of $\ve$ can be
chosen and to prove the boundedness of the family of controls
$(v_\varepsilon)_{\varepsilon\in(0,1)}$ in $L^2(0,T)$, the
dependence in $\ve$ of the constants $C$ and $\beta$ from
\eqref{ein} is required. This represents one of the most difficult
tasks of our work. We shall prove that $C$ and $\beta$ from
\eqref{ein} can be chosen independent of $\ve$, fact that ensures
the uniform boundedness of the sequence
$(v_\varepsilon)_{\varepsilon\in(0,1)}$ and the possibility to
pass to the limit as $\ve$ tends to zero in \eqref{ec.in}. The
main result of this paper reads as follows.

\begin{teor}\label{theo-main} Let $\al\in[0,1)\setminus\left\{\frac{1}{2}\right\}$ and
$f\in L^2(0,\pi)$ be a function such that $\ds\widehat{f}_n\neq0$
for every $n\geq 1$. There exists a time $T>0$ with the property
that, for any $(u^0,u^1)\in {\cal H}_0$ and $\varepsilon\in(0,1)$,
there exists a control $v_\varepsilon\in L^2(0,T)$ of
\eqref{ec.in} such that the family $(v_\varepsilon)_{\varepsilon
\in(0,1)}$ is uniformly bounded in $L^2(0,T)$ and any weak limit
$v$ of it, as $\varepsilon$ tends to zero, is a control in time
$T$ for equation \eqref{ec.in0}.
\end{teor}

The controllability problem studied in this paper belongs to the
interface between parabolic and hyperbolic equations. From this
point of view, it is related to \cite{Cor,Gla,LM}, where the
controllability of the transport equation is addressed after the
introduction of a vanishing viscosity term. In \cite{Cor} Carleman
estimates are used to obtain an uniform bound for the family of
controls. The same result is shown in \cite{Gla}, improving the
control time, by means of nonharmonic Fourier analysis and
biorthogonal technique. The recent article \cite{LM} deals with a
nonlinear scalar conservation law perturbed by a small viscosity
term and proves the uniform boundedness of the boundary controls.
Related problems in which controls for an equation are obtained as
limits of controls of equations of different type may be also
found in \cite{Lop,MR,Mil,Rus,Zua}.

In order to justify the damping mechanism introduced in
\eqref{ec.in}, which involves the fractional power $\alpha$ of the
Laplace operator, let us point out that sometimes it may be useful
to control the amount of dissipation introduced in the system not
only by means of the vanishing parameter $\varepsilon$ but also by
an adequate choice of the differential operator. For instance, the
convergence rates in some perturbed problems can be improved by
choosing a viscosity operator of lower order (see, for instance,
\cite{CI} in the context of Hamilton-Jacobi equations). In
\eqref{ec.in} this is achieved through the parameter $\alpha$. The
case $\alpha=1$ has been studied, for a slightly different
problem, in \cite{MOP}, where a uniform controllability result
with respect to the viscosity is proved. Theorem \ref{theo-main}
shows that a similar result holds for any $\alpha\in
\left[0,1\right)\setminus\left\{\frac{1}{2}\right\}$. Note that,
if $\alpha\in \left[0,\frac{1}{2}\right)$, the imaginary parts of
the eigenvalues $\lambda_n$ dominate the real ones and problem
\eqref{ec.in} has the same hyperbolic character as in the limit
case $\ve=0$. On the contrary, if $\alpha\in
\left(\frac{1}{2},1\right)$, \eqref{ec.in} has a parabolic type.
In this case we are dealing with a truly singular control problem
and the pass to the limit is sensibly more difficult. Finally, let
us remark that $\alpha=\frac{1}{2}$ is a singular case in which
the basic controllability properties (such as spectral
controllability) of \eqref{ec.in} do not hold.

For $\alpha\in[0,1)\setminus\{\frac{1}{2}\}$, the construction
from the proof of Theorem \ref{theo-main} implies that the
following Ingham-type inequality (see \cite{Ing2}) holds, for any
finite sequence $(\beta_n)_{n\in\mathbb{Z}^*}$ and $T$
sufficiently large,
\begin{equation}\label{ineg1}
C(T,\alpha)\sum_{n\in\mathbb{Z}^*}|\beta_n|^2e^{-\omega\ve
|n|^{2\alpha}}\leq
\int_{-T}^T\left|\sum_{n\in\mathbb{Z}^*}\beta_ne^{\lambda_n
t}\right|^2dt,
\end{equation}
where $\ve\in (0,1)$, $\omega$ is an absolute positive constant
and $C$ a positive constant depending of $T$ and $\alpha$ but
independent of $\ve.$ From this point of view our article extends
the results from \cite{E,Han,S}, where Ingham-type inequalities
are obtained under a more restrictive uniform sparsity condition
of the sequence $(\lm_n)_{n\in\mathbb{Z}^*}$. Indeed, one of the
major difficulty in our study is related to the fact that the
sequence of our eigenvalues $(\lambda_n)_{n\in\mathbb{Z}^*}$  is
not included in a sector of the positive real axis and does not
verify a uniform separation condition of the type
$$|\lambda_n-\lambda_m|\geq \delta |n^\beta-m^\beta|\qquad (n,m\in\mathbb{Z}^*),$$
for some $\beta>1$ and $\delta>0$ independent of $\ve$. The fact
that $C(T,\alpha)$ in \eqref{ineg1} does not depend of $\ve$ is of
fundamental importance since it ensures the uniform boundedness of
a family of controls $(v_{\ve})_{\ve\in(0,1)}$ for \eqref{ec.in}
and the possibility to pass to the limit in order to obtain a
control $v$ for \eqref{ec.in0}.

The rest of the paper is organized as follows. Section 2 gives the
equivalent characterization of the controllability property in
terms of a moment problem. The core of the paper is Section 3
where two biorthogonal sequences to the family $\left(e^{\lm_n
t}\right)_{n\in\mathbb{Z}^*}$ are constructed and evaluated. The
proof of Theorem \ref{theo-main} is provided in Section 4. The
article ends with an Appendix in which a technical lemma is
proved.

\section{The moment problem}\label{sec2}
In this section we show the equivalence between the
controllability problem \eqref{ec.in} $-$ \eqref{ec.coco} and the
moment problem \eqref{pm}. In order to do this we need first a
result concerning the existence of solutions for equation
\eqref{ec.in}. More precisely we have the following property.

\begin{prop}\label{th1}
Given any $T>0$, $\ve\in(0,1)$, $h\in L^{2}(0,T;L^{2}(0,\pi))$ and
$(u^0,u^1)\in H_0^1(0,\pi)\times L^{2}(0,\pi)$, there exists a
unique weak solution $(u,u_t)\in C([0,T],H_0^1(0,\pi)\times
L^{2}(0,\pi))$ of the problem
\begin{equation}\label{ec2.1}
\left\{
\begin{array}{ll}
\ds u_{tt}+(-\partial^{2}_{xx})u+
\ve^2(-\partial^{2}_{xx})^{2\al}u+2\ve(-\partial^{2}_{xx})^{\al}u_t=h(t,x)&(x,t)\in(0,\pi)\times(0,T)\\
\ds u(t,0)=u(t,\pi)=0&t\in(0,T)\\
\ds u(0,x)=u^{0}(x)\qquad u_t(0,x)=u^1(x)&x\in(0,\pi).
\end{array}\right.
\end{equation}
\end{prop}

\demo Since the operator $(D(A),A),$ where
$D(A)=\left\{\begin{array}{lcl}D((-\partial^{2}_{xx})^{\al})&\text{if}&\ve>0\\\\
D(-\partial^{2}_{xx})&\text{if}&\ve=0
\end{array}\right.$ and

$A=\left(\begin{array}{cc}0 & -I
\\-\partial_{xx}^2 +\varepsilon^2 (-\partial_{xx}^2)^{2\alpha}&
2\ve(-\partial_{xx}^2)^\alpha\end{array}\right),$ is maximal and monotone in $H_0^1(0,\pi)\times L^{2}(0,\pi),$ we apply the classical theory for the semigroups of contraction in $H^1_0(0,\pi)\times L^2(0,\pi)$ (see \cite{CH}).\ctd

Now we can give the  characterization of the  controllability
property of \eqref{ec.in} $-$ \eqref{ec.coco} in terms of a moment
problem. Based on the Fourier expansion of solutions, the moment
problems have been widely used in linear control theory. We refer
to \cite{Avd, Kom, Tuc, Zab} for a detailed  discussion of the
subject.

\begin{teor}\label{th2}
Let $T>0,\;\ve\in(0,1),\;(u^0,u^1)\in {\cal H}_0\;\text{and}\;f\in
L^{2}(0,\pi)$. There exists a control $v_{\ve}\in L^{2}(0,T)$ such
that the solution $(u,u_t)$ of equation (\ref{ec.in}) verifies
\eqref{ec.coco}, if and only if, $v_{\ve}\in L^{2}(0,T)$ satisfies
\begin{equation}\label{ec2.2}
\widehat{f}_{|n|}\int_{-\frac{T}{2}}^{\frac{T}{2}}v_{\ve}\left(t+\frac{T}{2}\right)
e^{\overline{\lambda}_{n}t}dt=
 -e^{-\overline{\lm}_{n}\frac{T}{2}}
 \left(\widehat{u}^1_{|n|}+\lm_n\widehat{u}^0_{|n|}\right)
 \qquad(n\in\mathbb{Z}^*),
\end{equation}
where $\lambda_{n}=in+\ve |n|^{2\al},$ for any $n\in\mathbb{Z}^*.$
\end{teor}
\demo We consider the ``adjoint" equation
\begin{equation}\label{ec2.3}
\left\{
\begin{array}{ll}
\ds \varphi_{tt}+(-\partial^{2}_{xx})\varphi+
\ve^2(-\partial^{2}_{xx})^{2\al}\varphi-2\ve(-\partial^{2}_{xx})^{\al}\varphi_t=0&(x,t)\in(0,\pi)\times(0,T)\\
\ds \varphi(t,0)=\varphi(t,\pi)=0&t\in(0,T)\\
\ds \varphi(T,x)=\varphi^{0}(x)\qquad
\varphi_t(T,x)=\varphi^1(x)&x\in(0,\pi).
\end{array}\right.
\end{equation}
If we multiply \eqref{ec.in}  by  $\overline{\varphi}$ and we
integrate by parts over $(0,T)\times(0,\pi)$,we deduce that
$v_{\ve}\in L^{2}(0,T)$ is a control for \eqref{ec.in} if, and
only if, it verifies
\begin{equation}\label{ec2.4}
\int_{0}^{T}v_{\ve}(t)\int_{0}^{\pi}f(x)
\overline{\varphi}(t,x)dxdt=-\int_{0}^{\pi}u^{1}(x)\overline{\varphi}(0,x)dx+
\int_{0}^{\pi}u^0(x)\left(\overline{\varphi}_t(0,x)-2\ve
\left(-\partial_{xx}^2\right)^{\al}\overline{\varphi}(0,x)\right)dx,
\end{equation}
for every $(\varphi,\varphi_t)$ solution of (\ref{ec2.3}) with the
initial data $(\varphi^0,\varphi^1)$. Since
$\left(\sin(nx)\right)_{n\geq1}$ is a basis for $L^2(0,\pi)$ we
have to check \eqref{ec2.4} only for the initial data of the form
$(\varphi^0,\varphi^1)=(\sin(nx),0)$ and
$(\varphi^0,\varphi^1)=(0,\sin(nx))$, for each $n\geq1$. In the
first case the solution of \eqref{ec2.4} is given by
\begin{equation}\label{unu}
\varphi(t,x)=
\left(\frac{\overline{\lm}_n}{\overline{\lm}_n-\lm_n}
e^{(t-T)\lm_n}+\frac{\lm_n}{\lm_n-\overline{\lm}_n}
e^{(t-T)\overline{\lm}_n}\right)\sin(nx)\qquad(n\in\mathbb{N}^*),
\end{equation}
whereas in the second case it becomes
\begin{equation}\label{doi}
\varphi(t,x)= \left(\frac{1}{\lm_n-\overline{\lm}_n}
e^{(t-T)\lm_n}+\frac{1}{\overline{\lm}_n-\lm_n}
e^{(t-T)\overline{\lm}_n}\right)\sin(nx)\qquad(n\in\mathbb{N}^*).
\end{equation}
By tacking in \eqref{ec2.4} $\varphi$ of the form \eqref{unu} and
\eqref{doi}, we obtained  that $v_{\ve}\in L^{2}(0,T)$ is a control
of \eqref{ec.in} if and only if it verifies \eqref{ec2.2}.\ctd

\begin{obs}
Note that $(\lm_n)_{n\in\mathbb{Z}^*}$ introduced in the previous
theorem are the eigenvalues of the differential operator
corresponding to the ``adjoint" equation \eqref{ec2.3}.
\end{obs}
\begin{obs}
The condition $\widehat{f}_n\neq0$ for any $n\in\mathbb{N}^*$ is
necessary in order to solve the moment problem \eqref{ec2.2} for
any initial data in ${\cal H}_0$. Indeed, if there exists
$n_0\in\mathbb{N}^*$ such that $\widehat{f}_{n_0}=0,$ then
\eqref{ec2.2} has a solution only if the initial data $(u^0,u^1)$
verify the additional condition
$\widehat{u}^1_{n_0}+\lm_{n_0}\widehat{u}^0_{n_0}=0$.
\end{obs}

We recall that  $(\theta_{m})_{m\in\mathbb{Z}^*}\in
L^{2}(-\frac{T}{2},\frac{T}{2})$ is a biorthogonal sequence to the
family of exponential functions
$\left(e^{\lambda_{n}t}\right)_{n\in\mathbb{Z}^*}\in
L^{2}(-\frac{T}{2},\frac{T}{2})$ if and only if $$
\int_{-\frac{T}{2}}^{\frac{T}{2}}\theta_{m}(t)
e^{\overline{\lambda}_{n}t}dt=\delta_{mn}\qquad(m,n\in\mathbb{Z}^*).$$

It is easy to see from \eqref{ec2.2} that,  if
$(\theta_{m})_{m\in\mathbb{Z}^*}$ is a biorthogonal sequence to
the family of exponential functions
$\left(e^{\lambda_{n}t}\right)_{n\in\mathbb{Z}^*}$ in
$L^{2}(-\frac{T}{2},\frac{T}{2}),$ then a control of \eqref{ec.in}
is given by
\begin{equation}\label{ec2.5}
v_\ve(t)=-\sum_{m\in\mathbb{Z}^*}\frac{e^{-\overline{\lm}_{m}\frac{T}{2}}}
{\widehat{f}_{|m|}}\left(\widehat{u}^1_{|m|}+\lm_m\widehat{u}^0_{|m|}\right)
\theta_{m}\left(t-\frac{T}{2}\right) \qquad (t\in(0,T)),
\end{equation}
provided that the right hand side converges in $L^{2}(0,T).$  Now
the main problem is to show that there exists a biorthogonal
sequence $(\theta_{m})_{m\in\mathbb{Z}^*}$ to the family of
exponential functions
$\left(e^{\lambda_{n}t}\right)_{n\in\mathbb{Z}^*}$ in
$L^{2}(-\frac{T}{2},\frac{T}{2})$ and to evaluate its norm, in
order to prove the convergence of the right hand side of
\eqref{ec2.5} for any $(u^0,u^1)\in{\cal H}_0$.

\section{Construction of a biorthogonal sequence}\label{sec3}
The aim of this section is to construct and evaluate an explicit
biorthogonal sequence to the family
$\left(e^{t\lambda_n}\right)_{n\in\mathbb{Z}^*}$ in
$L^2\left(-\frac{T}{2},\frac{T}{2}\right)$, where
$\lambda_{n}=in+\ve |n|^{2\al}$ are the eigenvalues introduced in
Theorem \ref{th2}. In order to do that, we define a family
$\left(\Psi_m(z)\right)_{m\in\mathbb{Z}^*}$  of entire functions
of exponential type independent of $\ve$ (see, for instance,
\cite{You}) such that $\Psi_m(i\overline{\lm}_n)=\delta_{mn}.$ The
inverse Fourier transform of
$\left(\Psi_m\right)_{m\in\mathbb{Z}^*}$ will give us the
biorthogonal sequence $(\theta_m)_{m\in\mathbb{Z}^*}$ that we are
looking for. Each $\Psi_m$ is obtained from a Weierstrass product
$P_m$ multiplied by an appropriate function $M_m$ with rapid decay
on the real axis. Such a method was used for the first time by
Paley and Wiener \cite{Wie} and, in the context of control
problems, by Fattorini and Russell \cite{Fat1,Fat2}. The main
difficulty in our analysis is to obtain good estimates for the
behavior of $P_m$ on the real axis and to construct an appropriate
multiplier $M_m$ in order to ensure the boundedness of $\Psi_m$ on
the real axis. As we shall see in Proposition \ref{lmprod} below,
the behavior of $\ln|P_m(x)|$ is always dominated by a subunitary
power of $|x|$, if $\al<\frac{1}{2}$ which facilitates the entire
construction and analysis. In the more difficult case
$\al>\frac{1}{2}$, $\ln|P_m(x)|$ behaves like $|x|$ on an interval
of length
$O\left(\left(\frac{1}{\ve}\right)^{\frac{1}{2\al-1}}\right)$. It
is precisely this property which makes the construction of $M_m$
more problematic and imposes the necessity of a careful analysis
of $P_m(x)$. Finally, the bounds obtained on the real axis for
$\Psi_m$ and the Plancherel's Theorem, will provide the desired
estimates for $\|\theta_m\|_{L^2(0,T)}$ and their dependence of
the parameters $m$, $\ve$ and $\alpha$.

\subsection{An entire function}\label{sec31}

In this subsection we construct the Weierstrass product $P_m$
mentioned above and we study some of its properties. For every
$m\in\mathbb{Z}^{*}$, we define the function
\begin{equation}\label{3.1}
P_{m}(z)=\prod\limits_{\substack{n\in\mathbb{Z}^{*}\\n \neq m}}
\left(1+\frac{zi}{\overline{\lambda}_{n}}\right)
\left(\frac{\overline{\lambda}_{n}}{\overline{\lambda}_{n}-\overline{\lambda}_{m}}\right).
\end{equation}

Firstly, let us state the following technical result concerning
the second part of the product $P_m$,  whose proof will be given
in the Appendix.

\begin{lema}\label{majqm}
There exists a constant $C>0$ such that, for all $\ve\in(0,1)$ and
$m\in\mathbb{Z}^{*}$, we have
\begin{equation}\label{eq00}
\prod\limits_{\substack{n\in\mathbb{Z}^{*}\\n\neq m}}
\left|\frac{\lambda_{n}}{\lambda_{n}-\lambda_{m}}\right|\leq16\exp(C
\ve m^{2\al}).
\end{equation}
\end{lema}

Now we pass to study the basic properties of the product $P_m$.

\begin{prop}\label{lm1}
Let $\alpha\in[0,1)\setminus\left\{\frac{1}{2}\right\}$ and
$\ve\in(0,1)$. For each $m\in \mathbb{Z}^*$, $P_{m}$ is an entire
function of exponential type at most $L_1$, where
$$L_1:=\left\{ \begin{array}{ll}
\max\left\{\frac{\sqrt{2}\pi}{2},\frac{4\ve}{1-2\alpha}\right\}&
\al\in\left[0,\frac{1}{2}\right)\\
\\
\max\left\{\frac{\sqrt{2}\pi}{2},\frac{8}{2\alpha-1}\right\}&
\al\in\left(\frac{1}{2},1\right),\end{array}\right.$$ with the
property that
\begin{equation}\label{3.2}
P_{m}(i\overline{\lambda}_{n})=\delta_{mn}\qquad
(n\in\mathbb{Z}^{*}).
\end{equation}
\end{prop}
\begin{obs}
Note that Proposition \ref{lm1} does not consider  the case
$\alpha=\frac{1}{2}$. In fact, if $\alpha=\frac{1}{2}$, the family
of exponential functions
$\left(e^{\lm_nt}\right)_{n\in\mathbb{Z}^*}$ is complete in
$L^2(0,a)$, for any $a>0.$ Indeed, since
\begin{equation}
\sum_{n\in\mathbb{Z}^*}\frac{\Re(\lm_n)}{1+|\lm_n|^2}=\infty,
\end{equation}
the completeness is a consequence of the Theorem
Sz\'{a}sz-M\"{u}ntz \cite{Sza}. Since this property remains true
if we eliminate a finite number of elements, we deduce that
$\left(e^{\lm_nt}\right)_{n\in\mathbb{Z}^*}$ is not minimal in
$L^2(0,a)$ and there exists no biorthogonal sequence to it in
$L^2(0,a)$. From the controllability point of view, it follows
that \eqref{ec.in} is not spectrally controllable if
$\alpha=\frac{1}{2}$.
\end{obs}

{\em Proof of Proposition \ref{lm1}.} By taking into account the
estimate \eqref{eq00} from Lemma \ref{majqm}, we only have to
study the function
$$E_{m}(z)=\prod\limits_{\substack{n\in\mathbb{Z}^{*}\\ n\neq
m}} \left(1+\frac{zi}{\overline{\lambda}_{n}}\right).$$ We have
that

$
\begin{array}{lcl}
\ds |E_{m}(z)|&=&\ds
\left|1+\frac{zi}{\lambda_{m}}\right|\prod\limits_{\substack{n\in\mathbb{N}^{*}\\n\neq
|m|}}
\left|\left(1+\frac{zi}{\overline{\lambda}_{n}}\right)\left(1+\frac{zi}{\lambda_{n}}\right)\right|
\leq \left(1+\frac{|z|}{|\lambda_{m}|}\right) \prod_{n=1}^{\infty}
\frac{|\lambda_{n}|^{2}+2|\Re(\lambda_{n})||z|+|z|^{2}}{|\lambda_{n}|^{2}}=\\
\\&=&\ds\left(1+\frac{|z|}{|\lambda_{m}|}\right)\exp\bigg(\underbrace{\sum_{n=1}^{[N_z]}\ln\left(1+\frac{2|\Re(\lambda_{n})||z|+|z|^{2}}
{|\lambda_{n}|^{2}}\right)}_{S_1}+\underbrace{\sum_{n=[N_z]+1}^{\infty}\ln\left(1+\frac{2|\Re(\lambda_{n})||z|+|z|^{2}}
{|\lambda_{n}|^{2}}\right)}_{S_2}\bigg),
\end{array}
$

where $N_z=\left(\frac{|z|}{2\ve}\right)^{\frac{1}{2\al}}.$ It follows that

$
\begin{array}{lll}
S_{1}&\leq&\ds\sum_{n=1}^{\left[N_z\right]}
\ln\left(1+\frac{2|z|^{2}}{|\lambda_{n}|^{2}}\right)\leq
\int_{0}^{
N_z}\ln\left(1+\frac{2|z|^{2}}{t^{2}+\ve^{2}t^{4\al}}\right)dt\leq
\int_{0}^{ N_z}\ln\left(1+\frac{2|z|^{2}}{t^{2}}\right)dt \leq\\
\\
&\leq &\sqrt{2}|z| \ds \int_0^\infty
\ln\left(1+\frac{1}{t^{2}}\right)dt=\frac{\sqrt{2}\pi}{2}|z|.
\end{array}
$

Thus, we have that
\begin{equation}\label{3.3}
S_{1}\leq \frac{\sqrt{2}\pi}{2} |z|.
\end{equation}
For $\al\in\left[0,\frac{1}{2}\right)$ we have that
$$
S_{2}\leq\ds
\ds4|z|\sum_{n=\left[ N_z\right]+1}^{\infty}
\frac{\ve n^{2\al}}{n^{2}+\ve^{2}n^{4\al}}\leq
4|z|\sum_{n=2}^{\infty}
\ve n^{2\al-2}\leq\frac{4\ve}{1-2\al}|z|.
$$

For $\al\in\left(\frac{1}{2},1\right)$  we define
$\ds\gv=\left(\frac{1}{\ve}\right)^{\frac{1}{2\al-1}}$ and we
deduce that

$
\begin{array}{lcl}
S_{2}&\leq&\ds
\sum_{n=\left[ N_z\right]+1}^{\infty}
\frac{4|\Re(\lambda_{n})||z|}{|\lambda_{n}|^{2}}=4|z|\sum_{n=\left[ N_z\right]+1}^{\infty}
\frac{\ve n^{2\al}}{n^{2}+\ve^{2}n^{4\al}}\leq4|z|\left(\sum_{n=1}^{[\gv]}+\sum_{n=[\gv]+1}^{\infty}
\right)\frac{\ve n^{2\al}}{n^{2}+\ve^{2}n^{4\al}}\leq\\
\\
&\leq&\ds
4|z|\left(\sum_{n=1}^{[\gv]}\frac{\ve}{n^{2-2\al}}+
\sum_{n=[\gv]+1}^{\infty}\frac{1}{\ve n^{2\al}}\right)\leq\frac{8}{2\al-1}|z|.
\end{array}
$

It follows that
\begin{equation}\label{3.4}
S_{2}\leq\left\{\begin{array}{ll}
\frac{4\ve}{1-2\al}|z|&\al\in\left[0,\frac{1}{2}\right)\\
\\
\frac{8}{2\al-1}|z|&\al\in\left(\frac{1}{2},1\right).\end{array}\right.
\end{equation}

From \eqref{3.3} and \eqref{3.4} we deduce that $E_m(z)$ is an
entire function of exponential type at most $L_1$ and the proof
ends.\ctd

\subsection{Evaluation of $P_m$ on the real axis}

This subsection is devoted to study the behavior of the entire
function $P_m$ on the real axis. The main result will be presented
in Proposition \ref{lmprod} below. Let us begin with the following
two simple lemmas. We recall that, for $\alpha>\frac{1}{2}$, we
have introduced the notation
$\ds\gv=\left(\frac{1}{\ve}\right)^{\frac{1}{2\al-1}}$.

\begin{lema}\label{lm5}
Let  $\ve\in(0,1)$ be fixed and
$\al\in\left[0,1\right)\setminus\left\{\frac{1}{2}\right\}$. For
any $x\geq0$ there exists a unique $\xe\geq0$ such that
$x^{2}=\xe^{2}+\ve^{2}\xe^{4\al}.$ Moreover, if
$\al\in\left[0,\frac{1}{2}\right)$ or
$\al\in\left(\frac{1}{2},1\right)$ and $x\leq\gv,$ then
\begin{equation}\label{majx}
\xe\leq x\leq\sqrt{2}\xe,
\end{equation}
\begin{equation}\label{ip1}
\left|x-|\lm_n|\right|\geq\frac{|\xe-n|}{\sqrt{2}}\qquad(n\in\mathbb{N}^*,\;n\leq\gv).
\end{equation}
Finally, if $\al\in\left(\frac{1}{2},1\right)$ and $x> \gv,$ then
\begin{equation}\label{ip2}
|\lm_n|-x\geq\frac{\ve(n^{2\al}-\xe^{2\al})}{2\sqrt{2}}\qquad(n\in\mathbb{N}^*,\;n>\gv).
\end{equation}
\end{lema}
\demo Let us first note that, for all $t\geq0,$ the
equation
\begin{equation}\label{3.7}
r^{2}+\ve^{2}r^{4\al}=t
\end{equation}
has only one solution in $[0,\infty)$. Indeed, if we define the function $f:[0,\infty)\to\mathbb{R}$,
$f(r)=r^{2}+\ve^{2}r^{4\al}-t$, it results that $f$ is increasing. Therefore
equation (\ref{3.7}) has  at most one solution in
$[0,\infty)$. On the other hand, we notice that
$f(0)=-t\;\text{and}\;\ds\lim_{r\ra\infty}f(r)=\infty,$ from which
we conclude that equation \eqref{3.7} has a unique solution in $[0,\infty)$. Concerning
\eqref{majx}, it is  obviously that $\xe\leq x$ and for the second
part of the inequality we notice that
$$x^2=\xe^2+\ve^2\xe^{4\al}\leq2\xe^2$$
for any $\al\in\left[0,\frac{1}{2}\right)$ or
$\al\in\left(\frac{1}{2},1\right)$ and $\xe\leq\gv.$ Finally,
taking into account that
$$|x-|\lm_n||=\left|\frac{\xe^2-n^2+\ve\left(\xe^{4\al}-n^{4\al}\right)}
{\sqrt{\xe^2+\ve^2\xe^{4\al}}+\sqrt{n^2+\ve^2n^{4\al}}}\right|,$$
relations \eqref{ip1} and \eqref{ip2} follows immediately. \ctd

\begin{lema}\label{majfrac}
The following inequalities hold
\begin{equation}\label{majrap}
\frac{n^{4\alpha-2}-x^{4\alpha-2}} {n^{2}-x^{2}}\leq \left\{
\begin{array}{ll}
 x^{4\alpha-4} & n\leq x\\
\\
n^{4\alpha-4} & 0\leq x\leq n
\end{array}
\right. \qquad
\begin{array}{l}\left(\alpha\in\left(\frac{1}{2},1\right)\right),\end{array}
\end{equation}

\begin{equation}
\frac{n^{4\alpha-2}-x^{4\alpha-2}} {n^{2}-x^{2}}\leq0 \qquad
\begin{array}{l}\left(\al\in\left[0,\frac{1}{2}\right)\right).\end{array}
\end{equation}
\end{lema}

\demo We notice that, when $x\geq n,$ then
$$\ds\frac{n^{4\al-2}-x^{4\al-2}}{n^{2}-x^{2}}=x^{4\al-4}
\frac{\left(\frac{n}{x}\right)^{4\al-2}-1}{\left(\frac{n}{x}\right)^{2}-1}\leq
\left\{
\begin{array}{ll}
0&\al\in\left[0,\frac{1}{2}\right)\\
\\
x^{4\al-4}&\al\in\left(\frac{1}{2},1\right),
\end{array}
\right.$$ The case $n\leq x$ is treated similarly.\ctd

\

The main result from this subsection is the following estimate  of
the function $P_m$ on the real axis.

\begin{prop}\label{lmprod}
Let $\ve\in(0,1)$ and $m\in\mathbb{Z}^*$. For each
$\al\in\left[0,1\right)\setminus\frac{1}{2}$ there exist two
positive constants $C$ and $\omega,$ independent of $\ve$ and $m,$
such that the function $P_m$ defined by \eqref{3.1} verifies
\begin{equation}\label{3.5}
|P_{m}(x)|\leq C\exp\left[\, \omega\left(
\varphi_{\ve}(x)+|\Re
(\lm_m)|\right)\right]\qquad(x\in\mathbb{R}),
\end{equation}
where
\begin{equation}\label{defphi}
\varphi_{\ve}(x)=\left\{
\begin{array}{lcl}
\ve|x|^{2\al}&\,\text{if}\,&\left(\al\in\left[0,\frac{1}{2}\right)\text{ and }
x\in\mathbb{R}\right)\quad or \quad\left(\al\in
\left(\frac{1}{2},1\right)\text{ and }|x|\leq\gv\right)\\
\\
\left(\frac{|x|}{\ve}\right)^{\frac{1}{2\al}}&\,\text{if}\,&\al\in\left(\frac{1}{2},1\right)\text{ and }
|x|>\gv.
\end{array}
\right.
\end{equation}
\end{prop}

\demo  With the notations from Proposition \ref{lm1} and by taking
into account estimate \eqref{eq00} from Lemma \ref{majqm}, it follows that it is
enough to evaluate $|E_{m}(x)|.$ Moreover, since $P_m$ is a
continuous function it is sufficient to consider  $x\neq
|\lambda_n|$ for all $n\in\mathbb{Z}^*$. In the sequel, $C$
denotes a generic constant which may change from one row to
another but it is always independent of $\ve$ and $m.$

To begin with, we evaluate $E_{m}$ on the real axis in the case
 $\left(\al\in\left[0,\frac{1}{2}\right)\text{ and } x\in\mathbb{R}\right)$
 or $\left(\al\in\left(\frac{1}{2},1\right)\text{ and }|x|\leq\gv\right).$

$
\begin{array}{lcl}
|E_{m}(x)|^{2}&=&\ds\left|\prod\limits_{\substack{n\in\mathbb{Z}^{*}\\n\neq
m}}\left(1+\frac{xi}{\overline{\lambda}_{n}}\right)\right|^{2}=
\left|1+\frac{xi}{\lm_m}\right|^2
\ds\prod\limits_{\substack{n=1\\n\neq m}}^{\infty}\left|\frac{|\lambda_{n}|^{2}+2xi\Re (\lambda_{n})-x^{2}} {|\lambda_{n}|^{2}}\right|^{2}=\\
\\
&=&\ds\frac{(x+m)^2+\ve^2m^{4\al}}{m^2+\ve^2m^{4\al}}
\underbrace{\prod\limits_{\substack{n=1\\n\neq m}}^{\infty}
\frac{\left(|\lambda_{n}|^{2}-x^{2}\right)^{2}}{|\lambda_{n}|^{4}}}_{E_{m}^{1}(x)}
\;\,\underbrace{\prod\limits_{\substack{n=1\\n\neq m}}^{\infty}
\frac{\left(|\lambda_{n}|^{2}-x^{2}\right)^{2}+4x^{2}\left(\Re
\lambda_{n}\right)^{2}}
{\left(|\lambda_{n}|^{2}-x^{2}\right)^{2}}}_{E_{m}^{2}(x)}.
\end{array}
$

We shall consider that $x\geq0.$ The opposite case can be treated
in a similar way. Now, we evaluate $E_{m}^{1}(x)$ by using Lemma
\ref{lm5}. We have that
$$
\begin{array}{lcl}
\left|E_{m}^{1}(x)\right|&=&\ds\prod\limits_{\substack{n=1\\n\neq
m}}^{\infty}
\left|\frac{|\lambda_{n}|^{2}-x^{2}}{|\lambda_{n}|^{2}}\right|^{2}=
\prod\limits_{\substack{n=1\\n\neq m}}^{\infty}
\left|\frac{n^{2}+\ve^{2}n^{4\al}-\xe^{2}-\ve^{2}\xe^{4\al}}
{n^{2}+\ve^{2}n^{4\al}}\right|^{2}=\\
\\
&=&\ds\left(\frac{m^{2}}{m^{2}-\xe^{2}}\right)^{2}
\prod_{n=1}^{\infty}\left|\frac{n^{2}-\xe^{2}}{n^{2}}\right|^{2}
\underbrace{\prod\limits_{\substack{n=1\\n\neq m}}^{\infty}
\left|1+\frac{\ve^{2}\xe^{2}}{1+\ve^{2}n^{4\al-2}}
\frac{n^{4\al-2}-\xe^{4\al-2}}{n^{2}-\xe^{2}}\right|^{2}}_{A_m^1(x)}=\\
\\
&=&\ds\left(\frac{m^{2}}{m^{2}-\xe^{2}}\right)^{2}\left(\frac{\sin\pi\xe}{\pi\xe}\right)^{2}
A_{m}^1(x)\leq\frac{m^{2}}{(m-\xe)^{2}}\left(\frac{\sin\pi\xe}{\pi\xe}\right)^{2}
A_{m}^1(x).
\end{array}
$$

From Lemma \ref{majfrac} we deduce that, if
$\al\in\left[0,\frac{1}{2}\right)$, the product $A_m^1(x)$ has the
following property
$$A_m^1(x)\leq1.$$ On the other hand, if $\al\in\left(\frac{1}{2},1\right)$, from
\eqref{majrap} we deduce that

$
\begin{array}{lcl}
A_{m}^1(x)&\leq&\ds\exp\left(2\sum_{n=1}^{[\xe]}\ln\left(1+\frac{\ve^{2}\xe^{4\al-2}}{1+\ve^{2}n^{4\al-2}}\right)+
2\sum_{n=[\xe]+1}^{\infty}\ln\left(1+\frac{\xe^{2}\ve^{2}n^{4\al-4}}{1+\ve^{2}n^{4\al-2}}\right)\right)\leq\\
\\
&\leq&\ds\exp\bigg(\underbrace{2\sum_{n=1}^{[\xe]}\frac{\ve^{2}\xe^{4\al-2}}{1+\ve^{2}n^{4\al-2}}}_{A_m^{11}(x)}+
\underbrace{2\sum_{n=[\xe]+1}^{\infty}\frac{\xe^{2}\ve^{2}n^{4\al-4}}{1+\ve^{2}n^{4\al-2}}}_{A_m^{12}(x)}\bigg).
\end{array}
$

Next, we proceed to evaluate the sums $A_m^{11}(x)$ and $A_m^{12}(x)$. Firstly, we have that
$$
A_m^{11}(x)=2\sum_{n=1}^{[\xe]}\frac{\ve^{2}\xe^{4\al-2}}{1+\ve^{2}n^{4\al-2}}\leq
2\int_{0}^{\xe}\frac{\ve^{2}\xe^{4\al-2}}
{1+\ve^{2}t^{4\al-2}}dt\leq2\int_{0}^{\xe}\ve^{2}\xe^{4\al-2}dt
\leq2\ve^{2}\xe^{4\al-1}\leq2\ve\xe^{2\al},$$and secondly  we deduce that
$$
A_m^{12}(x)=\ds\sum_{n=[\xe]+1}^{\infty}\frac{2\xe^{2}\ve^{2}n^{4\al-4}}{1+\ve^{2}n^{4\al-2}}
 \leq\int_{\xe}^{\infty}\frac{2\xe^{2}\ve^{2}t^{4\al-4}}{1+\ve^{2}t^{4\al-2}}dt=\ds
\int_{\xe}^{\gv}\frac{2\xe^{2}\ve^{2}t^{4\al-4}}{1+\ve^{2}t^{4\al-2}}dt+
\int_{\gv}^{\infty}\frac{2\xe^{2}\ve^{2}t^{4\al-4}}{1+\ve^{2}t^{4\al-2}}dt\leq
2\ve\xe^{2\al}.
$$
Consequently, we have proved that
\begin{equation}\label{majE1}
E_{m}^{1}(x)\leq
\left\{\begin{array}{lc}
\ds\frac{m^{2}}{(m-\xe)^{2}}\left(\frac{\sin\pi\xe}{\pi\xe}\right)^{2}&
\al\in\left[0,\frac{1}{2}\right),\;x\in\mathbb{R}\\
\\
\ds\frac{m^{2}}{(m-\xe)^{2}}\left(\frac{\sin\pi\xe}{\pi\xe}\right)^{2}
\exp\left(4\ve x^{2\al}\right)&
\al\in\left(\frac{1}{2},1\right),\;x\leq\gv.
\end{array}
\right.
\end{equation}

Nextly, we evaluate the product $E_{m}^2(x)$. In the following
estimates we shall use the notation $n_x=[\xe]$ and relations
\eqref{majx} $-$ \eqref{ip1} from Lemma \ref{lm5}.
$$
\begin{array}{c}
E_{m}^{2}(x)=\ds\underbrace{\frac{\left(|\lambda_{m}|^{2}-x^{2}\right)^{2}}
{\left(|\lambda_{m}|^{2}-x^{2}\right)^{2}+4x^{2}\Re^{2}(
\lambda_{m})} \left(1+\frac{4x^{2}\Re^{2}(
\lambda_{n_x})}{\left(|\lambda_{n_x}|^{2}-x^{2}\right)^{2}}\right)
\left(1+\frac{4x^{2}\Re^{2} (\lambda_{n_x+1})}
{\left(|\lambda_{n_x+1}|^{2}-x^{2}\right)^{2}}\right)}_{E_m^{21}(x)}\\
\\
\ds\prod_{n=1}^{n_x-1} \left(1+\frac{4x^{2}\Re^{2}(
\lambda_{n})}{\left(|\lambda_{n}|^{2}-x^{2}\right)^{2}}\right)
\prod_{n=n_x+2}^{\infty}\left(1+ \frac{4x^{2}\Re^{2}(
\lambda_{n})}{\left(|\lambda_{n}|^{2}-x^{2}\right)^{2}}\right)\leq
\\
\\
\ds\leq E_m^{21}(x) \prod_{n=1}^{n_x-1} \left(1+\frac{4\Re^{2}(\lambda_{n})}
{\left(|\lambda_{n}|-x\right)^{2}}\right)
\prod_{n=n_x+2}^{\infty}\left(1+ \frac{4x^{2}\Re^{2}(\lambda_{n})}
{|\lambda_n|^2\left(|\lambda_{n}|-x\right)^{2}}\right)=\\
\\
\ds=E_m^{21}(x)\exp\bigg[\underbrace{\sum_{n=1}^{n_x-1}\ln\left(1+\frac{4\Re^{2}(\lambda_{n})}
{\left(x-|\lambda_{n}|\right)^{2}}\right)}_{E_m^{22}(x)}+\underbrace{\sum_{n=n_x+2}^{\infty}\ln
\left(1+ \frac{4x^{2}\Re^{2}(\lambda_{n})}
{|\lambda_n|^2\left(|\lambda_{n}|-x\right)^{2}}\right)}_{E_m^{23}(x)}\bigg].
\end{array}
$$
We estimate $E_m^{22}(x)$ by using Lemma \ref{lm5} as follows
$$
E_m^{22}(x)\leq\sum_{n=1}^{n_x-1}\ln\left(1+\frac{8\ve^2n^{4\al}}
{\left(n_x-n\right)^{2}}\right)\leq32\al\int_0^{n_x}\frac{
n_x\ve^2t^{4\al-1}}{(n_x-t)^2+8\ve^2t^{4\al}}dt.$$ To bound from
above the last integral we have to split the interval $(0,n_x)$ in
three parts, by taking into account the following inequalities
\begin{equation}
\ds0\leq\frac{ n_x}{2}\leq n_x-\frac{\ve n_x^{2\al}} {1+\ve
n_x^{2\al-1}}\leq n_x.
\end{equation}
Thus, we have that
$$
\begin{array}{c}
\ds\int_0^{n_x}\frac{
n_x\ve^2t^{4\al-1}}{(n_x-t)^2+8\ve^2t^{4\al}}dt\leq\int_0^{\frac{n_x}{2}}\frac{\ve^2n_xt^{4\al-1}}{(n_x-t)^2}dt+
\int_{\frac{n_x}{2}}^{n_x-\frac{\ve n_x^{2\al}}{1+\ve
n_x^{2\al-1}}}\frac{\ve^2n_xt^{4\al-1}}{(n_x-t)^2}dt+
\int_{n_x-\frac{\ve n_x^{2\al}}{1+\ve n_x^{2\al-1}}}^{n_x}\frac{n_x}{t}dt\leq\\
\\
\ds\leq\int_0^{\frac{n_x}{2}}\frac{\ve^2n_xt^{4\al-1}}{\left(\frac{n_x}{2}\right)^2}dt+
\int_{\frac{n_x}{2}}^{n_x-\frac{\ve n_x^{2\al}}{1+\ve n_x^{2\al-1}}}
\frac{\ve^2n_x^{4\al}}{(n_x-t)^2}dt+
n_x\ln(1+\ve n_x^{2\al-1})\leq C\ve n_x^{2\al}.
\end{array}
$$
The last inequality takes place because $\ve^2n_x^{4\al-1}\leq\ve
n_x^{2\al}$ if $\left(\al\in\left[0,\frac{1}{2}\right)\text{ and }
x\in\mathbb{R}\right)$ or
$\left(\al\in\left(\frac{1}{2},1\right)\text{ and
}x\leq\gv\right).$ Now let us evaluate $E_m^{23}(x)$ by treating
separately the following cases:

Case I. $\al\in\left[0,\frac{1}{2}\right)$ and $x\in\mathbb{R}.$
Using relation \eqref{ip1} from Lemma \ref{lm5} we obtain that
$$
\begin{array}{c}
\ds E_m^{23}(x)\leq\int_{\xe+1}^{\infty} \ln\left(
1+\frac{8x^2\ve^2t^{4\al-2}}{(t-\xe)^2}\right)dt \leq
32\int_{\xe+1}^{\infty}\frac{x^2\ve^2t^{4\al-2}}{(t-\xe)^2+x^2\ve^2t^{4\al-2}}dt=\\
\\
\ds=32\left(\int_{\xe+1}^{2\xe}+\int_{2\xe}^{\infty}\right)
\frac{x^2\ve^2t^{4\al-2}}{(t-\xe)^2+x^2\ve^2t^{4\al-2}}dt\leq
C\ve\xe^{2\al}.
\end{array}
$$

Case II. $\al\in\left(\frac{1}{2},1\right)$ and $x\leq\gv.$ Using
the relations \eqref{ip1} and \eqref{ip2} from Lemma \ref{lm5} we
have that
$$
E_m^{23}(x)\leq\underbrace{\int_{\xe}^{\gv}\ln\left(1+\frac{8x^{2}\ve^{2}t^{4\al-2}}
{\left(t-\xe\right)^{2}}\right)dt}_{I_1}+
\underbrace{\int_{\gv}^{\infty}\ln\left(1+\frac{32x^{2}}
{\ve^{2}\left(t^{2\al}-\xe^{2\al}\right)^{2}}\right)dt}_{I_2} .
$$

To evaluate  $I_1$ we integrate by parts and we take into account that there exists a constant
$v\in(0,1)$ such that $\xe\leq v\gv$. The existence of this constant allows us to separate the interval
$(\xe,\gv)$ as follows
$$
\begin{array}{c}
I_1\leq\ds 4x^2\ve^\frac{1}{2\al-1}+
8x^2\ve^2\left(\int_{\xe}^{\frac{1}{v}\xe}+\int_{\frac{1}{v}\xe}^{\gv}\right)\frac{t^{4\al-2}}
{(t-\xe)^{2}+4x^2\ve^2t^{4\al-2}}dt\leq\\
\\
\leq \ds 4\ve x^{2\al}+\int_{\xe}^{\frac{1}{v}\xe}\frac{C\xe^{4\al}\ve^2}
{(t-\xe)^{2}+\xe^{4\al}\ve^2}dt+\int_{\frac{1}{v}\xe}^{\gv}\frac{C\xe^2\ve^2t^{4\al-4}}
{1+\xe^2\ve^2t^{4\al-4}}dt\leq\\
\\
\leq\ds 4\ve x^{2\al}+\int_{0}^{\infty}\frac{C\xe^{4\al}\ve^2}
{t^{2}+\xe^{4\al}\ve^2}dt+\int_{\frac{1}{v}\xe}^{\gv}C\xe^2\ve^2
t^{4\al-4}dt\leq C\ve x^{2\al}.
\end{array}$$
In order to estimate $I_2$ we remark that
$t^{2\al}-\xe^{2\al}\geq(1-v^{2\al})t^{2\al}.$ It follows that
$$
I_2\leq
\int_{\gv}^{\infty}\ln\left(1+\frac{32x^2}{(1-v^{2\al})^2\ve^2t^{4\al}}\right)dt
\leq\ds\int_{\gv}^{\infty}\frac{32x^2}{(1-v^{2\al})^2\ve^2t^{4\al}}dt\leq
\frac{32}{(4\al-1)(1-v^{2\al})^2}\ve x^{2\al}.
$$
Thus, we have that
\begin{equation}\label{majE2}
E_m^2(x)\leq E_m^{21}(x)\exp(C\ve
|x|^{2\al}).
\end{equation}
Note that, for any $m\in\mathbb{Z}^*,$ there exists a positive constant $\widetilde{C},$ independent of $m$ and $\ve,$ such that
\begin{equation}\label{defC}
\frac{(x+m)^2+\ve^2m^{4\al}}{m^2+\ve^2m^{4\al}}\frac{m^{2}}{(m-\xe)^{2}}\left(\frac{\sin\pi\xe}{\pi\xe}\right)^{2}
\left(1+\frac{4x^{2}\Re^{2}(\lambda_{n_x})}{\left(|\lambda_{
n_x}|^{2}-x^{2}\right)^{2}}\right) \left(1+\frac{4x^{2}\Re^{2}
(\lambda_{n_x+1})}{\left(|\lambda_{
n_x+1}|^{2}-x^{2}\right)^{2}}\right)\leq \widetilde{C}.
\end{equation}
Indeed, the terms $\frac{m^{2}}{(m-\xe)^{2}}$,
$1+\frac{4x^{2}\Re^{2}(\lambda_{n_x})}{\left(|\lambda_{
n_x}|^{2}-x^{2}\right)^{2}}$ and $1+\frac{4x^{2}\Re^{2}
(\lambda_{n_x+1})}{\left(|\lambda_{n_x+1}|^{2}-x^{2}\right)^{2}}$
explodes as $x$ tends to $m$, $\lm_{n_x}$ and $\lm_{n_x+1}$
respectively, but not simultaneously. This allows us to couple
them with the sine function
$\left(\frac{\sin\pi\xe}{\pi\xe}\right)^{2}$ in order to obtain a
bounded function. On the other hand, when $x$ tends to infinity we
couple the first two terms so that we obtain once again a bounded
function.

Consequently, from \eqref{majE1}, \eqref{majE2} and \eqref{defC}
it follows that, for every $\al\in [0,1)\setminus \{\frac{1}{2}\}$
there exists an absolute positive constant $C,$ such that for any
$\left(x\in\mathbb{R} \mbox{ if }\al\in
\left[0,\frac{1}{2}\right)\right)$ or $\left(|x|\leq\gv\text{ if
}\al\in \left(\frac{1}{2},1\right)\right)$, we have that
\begin{equation}\label{E_m1}
|E_m(x)|\leq\exp(C\ve x^{2\al}).
\end{equation}

To conclude the proof it remains to evaluate the product
$E_{m}(x)$ in the case $\al\in\left(\frac{1}{2},1\right)$ and
$|x|>\gv.$ Note that
$$
\begin{array}{c}
|E_{m}(x)|^{2}=\ds\left|\prod\limits_{\substack{n\in\mathbb{Z}^{*}\\n\neq
m}} \left(1+\frac{xi}{\overline{\lambda}_{n}}\right)\right|^{2}=
\left|1+\frac{xi}{\lm_m}\right|^2
\left|\prod\limits_{\substack{n=1\\n\neq m}}^{\infty}
\left(1+\frac{xi}{\lambda_{n}}\right)
\left(1+\frac{xi}{\overline{\lambda}_{n}}\right)\right|^{2}\leq\\
\\
\leq\ds\frac{(m+x)^2+\ve^{2}m^{4\al}}{m^2+\ve^2m^{4\al}}\exp\bigg[
\underbrace{\sum_{n=1}^{[\eta_\ve(x)]}\ln\left(1+\frac{x^{4}+4x^{2}\Re^{2}\left(
\lambda_{n}\right)} {|\lambda_{n}|^{4}}\right)}_{S_m^1(x)}+\underbrace{\sum_{n=[\eta_\ve(x)]+1}^{\infty}\ln\left(1+\frac{x^{4}+4x^{2}\Re^{2}\left(
\lambda_{n}\right)} {|\lambda_{n}|^{4}}\right)}_{S_m^2(x)}\bigg],
\end{array}
$$
where $\ds\eta_{\ve}(x)=\left(\frac{x}{\sqrt{2}\ve}\right)^{\frac{1}{2\al}}$.
Now, we evaluate the above sum as follows.

$
\begin{array}{lcl}
S_m^1(x)&\leq&\ds
\sum_{n=1}^{[\eta_{\ve}(x)]}\ln\left(1+\frac{x^{4}}{|\lambda_{n}|^{4}}\right)\leq
\int_{0}^{\eta_{\ve}(x)}\ln\left(1+\frac{x^{4}}
{\left(t^{2}+\ve^{2}t^{4\al}\right)^{2}}\right)dt=\int_{0}^{\gv}\ln\left(1+\frac{x^{4}}{t^{4}}\right)dt+\\
\\
&+&\ds
\int_{\gv}^{\eta_{\ve}(x)}\ln\left(1+\frac{x^{4}}{\ve^{4}t^{8\al}}\right)dt\leq
C\left(\frac{x}{\ve}\right)^{\frac{1}{2\al}},
\end{array}
$

$
\begin{array}{lcl}
S_m^2(x)&\leq&\ds\sum_{n=[\eta_{\ve}(x)]+1}^{\infty}\frac{x^{4}+4x^{2}
\Re^{2}\left(\lambda_{n}\right)}{|\lambda_{n}|^{4}}\leq
\sum_{n=[\eta_{\ve}(x)]+1}^{\infty}\frac{8x^{2}\Re^{2}\left(\lambda_{n}\right)}
{|\lambda_{n}|^{4}}=\sum_{n=[\eta_{\ve}(x)]+1}^{\infty}\frac{8x^{2}\ve^{2}n^{4\al}}
{\left(n^{2}+\ve^{2}n^{4\al}\right)^{2}}\leq\\
\\
&\leq&\ds\sum_{n=[\eta_{\ve}(x)]+1}^{\infty}
\frac{8x^{2}}{\ve^{2}n^{4\al}}
\leq\ds\int_{\eta_{\ve}(x)}^{\infty}\frac{8x^{2}}{\ve^{2}}t^{-4\al}dt\leq
C\left(\frac{x}{\ve}\right)^{\frac{1}{2\al}}.
\end{array}
$

Thus, for $\al\in\left(\frac{1}{2},1\right)$ and $|x|>\gv$, we
have proved that
\begin{equation}\label{E_m2}
|E_{m}(x)|\leq\exp\left(C\left(\frac{|x|}{\ve}\right)^{\frac{1}{2\al}}\right).
\end{equation}
Now, by taking into account \eqref{E_m1} and \eqref{E_m2} the
proof of the Proposition ends.\ctd

\subsection{A multiplier}

In this subsection we construct a function, called multiplier,
used to compensate the grow of the product $P_m$ on the real axis
given in Proposition \ref{lmprod}.

Let $\varphi:[0,\infty)\to [0,\infty)$ be a continuous, increasing
and onto function. We define the real sequence $(a_{n})_{n\geq 1}$
by
\begin{equation}\label{an} \varphi(e a_{n})=n\qquad (n\geq
1)\end{equation} and we suppose that the following properties
hold:
\begin{enumerate}

\item[(I1)] $\displaystyle \sum_{n\geq 1}\frac{1}{a_{n}}\leq
L_2<\infty$

\item[(I2)] $\displaystyle
 \sum_{n=n_{m}}^{\infty}\frac{1}{a_{n}^{2}}\leq D \frac{1+
|\Re(\lm_m)|}{|\lm_{m}|^{2}},$ \end{enumerate} where $L_2$ and $D$
are two positive constants and
\begin{equation}\label{nm}
n_{m}= \left[\varphi(e|\lm_{m}|)\right] + 1 \qquad (m\geq 1).
\end{equation}

We have the following result.

\begin{lema}\label{lm9}
Let $ x\geq a_{n_m}$ and define $n_{x}:=[\varphi(e x)]$. Then
\begin{equation}\label{fo1}
\sum_{j=n_{m}}^{n_{x}}\ln\left(\frac{a_{j}}{x}\right)=
-\int_{a_{n_{m}}}^{x}\frac{A(u)-n_{m}+1}{u}du,\end{equation} where
$A(u)=\#\{a_{n}\leq u\}=\left[\varphi_\varepsilon(eu)\right].$
\end{lema}

\demo We have that $$
\begin{array}{l}
\ds-\int_{a_{n_{m}}}^{x}\frac{A(u)}{u}du
=\ds-\sum_{j=n_{m}}^{n_{x}-1}\int_{a_{j}}^{a_{j+1}}\frac{A(u)}{u}du-
\int_{a_{n_{x}}}^{x}\frac{A(u)}{u}du=
\ds-\sum_{j=n_{m}}^{n_{x}-1}\int_{a_{j}}^{a_{j+1}}\frac{j}{u}du-
\int_{a_{n_{x}}}^{x}\frac{n_{x}}{u}du\\
\\
=\ds\ln\left(\prod_{j=n_{m}}^{n_{x}-1}\frac{a_{j}^{j}}{a_{j+1}^{j}}\frac{a_{n_{x}}^{n_{x}}}{x^{n_{x}}}\right)=
\ln\left(\frac{a_{n_{m}}^{n_{m}-1}}{x^{n_{m}-1}}\prod_{j=n_{m}}^{n_{x}}\frac{a_{j}}{x}\right)
\ds =-\int_{a_{n_{m}}}^{x}\frac{n_m-1}{u}du+
\sum_{j=n_{m}}^{n_{x}}\ln\left(\frac{a_{j}}{x}\right),
\end{array}
$$ which completes the proof of the lemma.\ctd

Now we can construct a multiplier function.

\begin{teor}\label{lm10} Let $\varphi:[0,\infty)\to [0,\infty)$ be a continuous,
increasing, onto function such that the sequence $(a_{n})_{n\geq
1}$ defined by \eqref{an} verifies (I1) and (I2). For each $m\geq
1$ there exists $M_{m}:\CC\ra\CC$ with the following properties:
\begin{enumerate}
\item $M_{m}$ is an entire function of exponential type $L_2$

\item $\ds\left|M_{m}(x)\right|\leq \exp\left(-\varphi
(|x|)+\varphi(e|\lm_{m}|)+1\right)$ for all $x\in\mathbb{R}$

\item $\ds\left|M_{m}(i\overline{\lambda}_{m})\right|\geq
\exp\left(-\frac{D}{6}\left(1+|\Re(\lambda_{m})|\right)\right)$.
\end{enumerate}
\end{teor}
\demo By adapting an idea from \cite{Ing1}, we define the function
$M_{m}:\CC\ra\CC$ as follows
\begin{equation}\label{3.23}
M_{m}(z)=\prod_{n=n_{m}}^{\infty}\frac{\sin\left(\frac{z}{a_{n}}\right)}{\frac{z}{a_{n}}},
\end{equation}
where the sequence $(a_n)_{n\geq 1}$ is given by \eqref{an} and
$n_m$ is defined in \eqref{nm}.

 $M_{m}$ is an entire function of exponential
type. Indeed, this is a consequence of property (I1) of the
sequence $(a_n)_{n\geq 1}$ and the following estimate which holds
for each $N>n_m$,
$$
\prod_{n=n_{m}}^{N}\left|
\frac{\sin\left(\frac{z}{a_{n}}\right)}{\frac{z}{a_{n}}}\right|\leq\prod_{n=n_{m}}^{N}
e^{\left|\frac{z}{a_{n}}\right|}=
e^{|z|\ds\sum_{n=n_{m}}^{N}\frac{1}{a_{n}}}\leq e^{L_2|z|}.
$$

To prove the second property of $M_m$ we need to analyze the
following two cases:
\begin{enumerate}
\item[] Case 1: $x\leq e a_{n_m}$. We deduce that $\varphi(x) \leq
\varphi(e a_{n_m})=n_m\leq \varphi(e|\lm_m|)+1$ and consequently
$$
\begin{array}{c}
\ds\left|M_{m}(x)\right|=\ds\prod_{n=n_{m}}^{\infty}\left|
\frac{\sin\left(\frac{x}{a_{n}}\right)}{\frac{x}{a_{n}}}\right|\leq1
\leq\exp\left(\varphi\left(e |\lm_{m}|\right)-\varphi(
x)+1\right).
\end{array}$$

\item[]Case 2: $e a_{n_m}< x$. From Lemma \ref{lm9} we deduce that
$$
|M_{m}(x)|=\prod_{n=n_m}^{\infty}\left|
\frac{\sin\left(\frac{x}{a_{n}}\right)}{\frac{x}{a_{n}}}\right|\leq
\prod_{n=n_m}^{n_{x}}\left|\frac{a_{n}}{x}\right|=
\exp\left(\sum_{n=n_{m}}^{n_{x}}\ln\left(\frac{a_{n}}{x}\right)\right)=
\exp\left(-\int_{a_{n_{m}}}^{x}\frac{A(u)-n_{m}+1}{u}du\right).
$$

Since $a_{n_m}<\frac{x}{e}$, it follows that
$$
|M_{m}(x)|\leq
\exp\left(-\int_{\frac{x}{e}}^{x}\frac{A(u)-n_{m}+1}{u}du\right)\leq
\exp\left(-\int_{\frac{x}{e}}^{x}\frac{\varphi(x)-1-n_{m}+1}{u}du\right)=
\exp\left(-\varphi(x) + n_{m}\right).
$$
Since $n_m=[\varphi (e|\lambda_m|)] +1\leq \varphi
(e|\lambda_m|)+1,$ the second property of $M_m$ is proved.

\end{enumerate}

To prove the third property of $M_m$, note that
$$a_{n_{m}}=\frac{1}{e}\varphi^{-1}(n_m)\geq \frac{1}{e}\varphi^{-1}(\varphi(e|\lm_m|))
=|\lm_{m}|$$ and consequently
$\left|\frac{\lm_{m}}{a_{n}}\right|\leq 1$ for each $n\geq n_m$.
It follows that
$$
\begin{array}{lll}
\ds|M_{m}(i\overline{\lm}_{m})|&=&\ds
\prod_{n=n_{m}}^{\infty}\left|\frac{\sin\left(\frac{i\overline{\lm}_{m}}{a_{n}}\right)}
{\frac{i\overline{\lm}_{m}}{a_{n}}}\right|\geq
\prod_{n=n_{m}}^{\infty}\frac{\sin\left|\frac{\overline{\lm}_{m}}{a_{n}}\right|}
{\left|\frac{\overline{\lm}_{m}}{a_{n}}\right|}\geq
\prod_{n=n_{m}}^{\infty}\left|1-\frac{1}{6}\frac{|\overline{\lm}_{m}|^{2}}{a_{n}^{2}}\right|=\\
\\
&=&\ds
\exp\left(\sum_{n=n_{m}}^{\infty}\ln\left(1-\frac{1}{6}\frac{|\lm_{m}|^{2}}{a_{n}^{2}}\right)
\right)
\geq\ds\exp\left(-\frac{|\lm_{m}|^{2}}{6}\sum_{n=n_{m}}^{\infty}\frac{1}{a_{n}^{2}}\right)
\end{array}
$$
By using property (I2) of the sequence $(a_n)_{n\geq 1}$, we
deduce that the third property of $M_m$ also holds and the proof
of the theorem ends.\ctd

\begin{prop}\label{lm20} For $\alpha\in\left(0,1\right)\setminus \left\{\frac{1}{2}\right\}$ and
$\varepsilon\in(0,1)$, let $\varphi_{\ve}:[0,\infty)\to [0,\infty)$ be
the function defined by \eqref{defphi}. For each $m\geq 1$ there exists
$M_{m}:\CC\ra\CC$ with the following properties:
\begin{enumerate}
\item $M_{m}$ is an entire function of exponential type $L_2$

\item $\ds\left|M_{m}(x)\right|\leq \exp\left(-\varphi_\ve
(x)+2e^2|\Re(\lm_m)|+1\right)$ for all $x\in\mathbb{R}$

\item $\ds\left|M_{m}(i\overline{\lambda}_{m})\right|\geq
\exp\left(-D\left(1+|\Re(\lm_m)|)\right)\right)$,
\end{enumerate}
where $L_2$ and $D$ are positive constants independent of $m$ and
$\ve.$
\end{prop}
\demo The existence of the function $M_m$ follows from Theorem
\ref{lm10} if we prove that the function $\varphi_{\ve}$ verifies
the hypothesis from Theorem \ref{lm10} and
\begin{equation}\label{majf}
\varphi_{\ve}(e|\lm_m|)\leq2e^2|\Re(\lm_m)|\qquad(m\geq1).
\end{equation}

From \eqref{defphi} we deduce immediately that
$\varphi_{\ve}:(0,\infty)\ra(0,\infty)$ is continues, increasing,
onto and \eqref{majf} is verified. Moreover, the sequence
$(a_n)_{n\geq1}$ defined by $a_{n}=\frac{1}{e}\varphi^{-1}(n)$
verifies the following properties.
\begin{itemize}

\item For $\al\in\left(0,\frac{1}{2}\right)$ we have that
$$\frac{1}{e}\sum_{n\geq1}\frac{1}{a_{n}}=\ve^{\frac{1}{2\al}}+
\ve^{\frac{1}{2\al}}\sum_{n\geq2}
\left(\frac{1}{n}\right)^{\frac{1}{2\al}}\leq\frac{4\al+1}{2\al}\ve^{\frac{1}{2\al}}.
$$
\item For $\al\in\left(\frac{1}{2},1\right)$ we have that
$$
\begin{array}{c}
\ds\frac{1}{e}\sum_{n\geq1}\frac{1}{a_{n}}=\sum_{n\geq1}\frac{1}{\vf(n)}\leq
\int_{0}^{\infty}\frac{1}{\vf(s)}ds=
\int_{0}^{\gv}\left(\frac{\ve}{s}\right)^{\frac{1}{2\al}}ds+
\int_{\gv}^{\infty}\frac{1}{\ve s^{2\al}}ds=\frac{2\al+1}{2\al-1}.
\end{array}
$$
\end{itemize}
By taking $L_2=\frac{4\al+1}{2\al}\ve^{\frac{1}{2\al}}e$ for
$\al\in\left(0,\frac{1}{2}\right)$ and $L_2=\frac{2\al+1}{2\al-1}e$ for
$\al\in\left(\frac{1}{2},1\right)$
it follows that hypothesis
$(I1)$ is verified.

On the other hand we have that
\begin{itemize}
\item For $\al\in\left(0,\frac{1}{2}\right)$ we have that
$$\frac{1}{e^2}\sum_{n=n_{m}}^{\infty}\frac{1}{a_{n}^{2}}
=\sum_{n=n_{m}}^{\infty}\left(\frac{\ve}{n}\right)^{\frac{1}{\al}}
\leq \int_{\varphi_{\ve}(e|\lm_m|)}^{\infty}
\left(\frac{\ve}{s}\right)^{\frac{1}{\al}}ds
\leq\ve|\lm_m|^{2\al-2}\leq
\frac{1}{|\lm_m|^{2}}\left(1+2^{\alpha}|\Re(\lm_m)|\right).$$

\item For $\al\in\left(\frac{1}{2},1\right)$ it follows that
$$\ds\frac{1}{e^2}\sum_{n=n_{m}}^{\infty}\frac{1}{a_{n}^{2}}=\sum_{n=n_{m}}^{\infty}\frac{1}{(\vf(n))^{2}}\leq
\frac{1}{(\vf(n_{m}))^{2}}+\int_{n_{m}}^{\infty}\frac{1}{(\vf(s))^{2}}ds.$$
Since $n_{m}\geq\varphi_{\ve}(e|\lm_{m}|)$ it results that
\begin{equation}\label{3.22}
\frac{1}{(\vf(n_{m}))^{2}}\leq\frac{1}{e^2|\lm_{m}|^{2}}\qquad
(m\in\mathbb{N}^{*}).
\end{equation}

On the other hand, if $e|\lm_m|\leq\gv$ we have that
$$
\begin{array}{lcl}
\ds
\int_{\varphi_{\ve}(e|\lm_m|)}^{\infty}\frac{1}{(\vf(s))^{2}}ds&=&
\ds\int_{\varphi_{\ve}(e|\lm_m|)}^{\gv}\left(\frac{\ve}{s}\right)^{\frac{1}{\al}}ds+
\int_{\gv}^{\infty}\frac{1}{\ve^{2}s^{4\al}}ds \leq \\ \\ &\leq &
\ds\frac{\al}{1-\al}\ve^{\frac{1}{\al}}\left(\ve
(e|\lm_m|)^{2\al}\right)^{\frac{\al-1}{\al}}+
\frac{1}{4\al-1}\frac{1}{\ve^{2}}\gv^{1-4\al} \leq \ds
\frac{4\al}{1-\al}\frac{1}{|\lm_{m}|^{2}}|\Re (\lambda_m)|
\end{array}
$$
and if $e|\lm_m|> \gv$
$$
\begin{array}{l}
\ds\int_{\varphi_{\ve}(e|\lm_m|)}^{\infty}\frac{1}{(\vf(s))^{2}}ds=\ds
\int_{\varphi_{\ve}(e|\lm_m|)}^{\infty}\frac{1}{\ve^{2}s^{4\al}}ds\leq\ds
\frac{1}{4\al-1}e^{\frac{1-4\al}{2\al}}\frac{1}{\ve^{2}}\left(\frac{|\lm_m|}{\ve}\right)^{\frac{1-4\al}{2\al}}
\leq\ds\frac{4e^2}{4\al-1}\frac{1}{|\lm_{m}|^{2}}|\Re
(\lambda_m)|.
\end{array}
$$
\end{itemize}
 By taking $D=2^{\al}e^2$ for
$\al\in\left(0,\frac{1}{2}\right)$ and $D=\frac{4\al}{1-\al}e^2$
for $\al\in\left(\frac{1}{2},1\right)$ it follows that hypothesis
$(I2)$ is verified and the proof of the proposition finishes.\ctd

\subsection{Two biorthogonal sequences}

Now we have all the ingredients needed to construct a biorthogonal
sequence $(\theta_m)_{m\in\mathbb{Z}^*}$ to the family
$\left(e^{\lambda_n t}\right)_{n\in\mathbb{Z}^*}$, by using the
method presented in the section's introduction.

\begin{teor}\label{bio1}
Let $\ve\in(0,1)$. There exist $\widetilde{T}>0$ independent of
$\ve$ and a biorthogonal sequence
$(\theta_{m})_{m\in\mathbb{Z}^{*}}$ to the family
$(e^{\lm_{n}t})_{n\in\mathbb{Z}^{*}}$ in $\lt$, with the following
property
\begin{equation}\label{3.17}
\|\theta_{m}\|_{\lt}\leq C \exp(\beta|\Re (\lm_{m})|)\qquad
(m\in\mathbb{Z}^{*}),
\end{equation}
where $C$ and $\beta$ are positive constants independent of  $m$
and $\ve.$
\end{teor}
\demo If $\alpha\neq 0$, for each $m\in\mathbb{Z}^{*},$ let $P_m$
and $M_{|m|}$ be the functions from Propositions \ref{lm1} and
\ref{lm20}, respectively. We define the function
\begin{equation}\label{t1}
\Psi_{m}(z)=P_{m}(z)\left(\frac{M_{|m|}(z)}{M_{|m|}(i\overline{\lm}_{m})}\right)^{\omega}
\frac{\sin(\delta(z-i\overline{\lm}_{m}))}{\delta(z-i\overline{\lm}_{m})},
\end{equation}
where $\delta>0$  is an arbitrary constant and $\omega$ is the
constant from Proposition \ref{lmprod}. Let
\begin{equation}
\theta_{m}(t)=\frac{1}{2\pi}\int_{\mathbb{R}}\Psi_{m}(x)e^{ixt}dx.
\end{equation}

From Propositions \ref{lm1} and \ref{lm20} we deduce that there
exists $\widetilde{T}= 2(L_1+\omega L_2 +\delta)$, independent of
$\ve$, such that $\Psi_{m}$ is an entire function of exponential
type
 $\frac{\widetilde{T}}{2}.$ Moreover, from the estimate of the function
  $P_{m}$ on the real axis given by
Proposition \ref{lmprod} and the properties of the function
$M_{|m|}$ from Proposition \ref{lm20}, we obtain that
$$
\begin{array}{lcl}
\ds\int_{\mathbb{R}}|\Psi_{m}(x)|^{2}dx&\leq&\ds
Ce^{2\omega(1+D)}\exp\left(2\omega( 1+2e^2+D) |\Re
(\lm_{m})|\right)
\int_{\mathbb{R}} \left|\frac{\sin(\delta(x-i\lm_{m}))}{\delta(x-i\lm_{m})}\right|^{2}dx\\
\\
&\leq &\ds\frac{C}{\delta}\exp\left( 2(\omega+2\omega e^2+\omega
D+\delta)|\Re (\lm_{m})|\right) \int_{\mathbb{R}}
\left|\frac{\sin\, t}{t}\right|^{2}dt\leq C \exp\left(2\beta|\Re
(\lm_{m})|\right),
\end{array}
$$
where $\beta$ is any number greater that $\omega+2\omega
e^2+\omega D+\delta.$

Now, by taking into account the properties of $\Psi_m$ and by
applying the Paley-Wiener Theorem, \cite[Theorem 18, Sec.
2.4]{You}, we deduce that $\theta_{m}\in\lt$. Moreover, from the
inverse Fourier transform property we obtained that
$(\theta_m)_{m\in\mathbb{Z}^*}$ is a biorthogonal sequence to
$\left(e^{\lm_nt}\right)_{n\in\mathbb{Z}^*}$ in $\lt$. Finally,
from Plancherel's Theorem we deduce that (\ref{3.17}) holds.

If $\alpha=0$, we take \begin{equation}\label{t2}
\Psi_{m}(z)=P_{m}(z)
\frac{\sin(\delta(z-i\overline{\lm}_{m}))}{\delta(z-i\overline{\lm}_{m})},
\end{equation}
where $\delta>0$  is an arbitrary constant. The same argument as
before allows us to end the proof of the theorem.\ctd

The following result gives the existence of a new biorthogonal
sequence with better  norm properties than the one from Theorem
\ref{bio1}.

\begin{teor}\label{theo-2} Let $\ve\in(0,1)$. There exist $T_0>0$ independent of $\ve$
and a biorthogonal sequence $(\zeta_m)_{m\in\mathbb{Z}^\ast}$ to
the family $\left(e^{\lambda_nt}\right)_{n\in\mathbb{Z}^\ast}$ in
$L^2\left(-\frac{T_0}{2},\frac{T_0}{2}\right)$, such that, for any
finite sequence $(c_m)_{m\in \mathbb{Z}^\ast}$, we have
\begin{equation}\label{est2}
\int_{-\frac{T_0}{2}}^{\frac{T_0}{2}}\left| \sum_{m\in\mathbb{Z}^\ast}
c_m \zeta_{m}(t)\right|^2 dt\leq C(T_0) \sum_{m\in \mathbb{Z}^\ast}
|c_m|^2 e^{2\beta |\Re(\lambda_m)|},
\end{equation}
where $\beta$ is the same as in Theorem \ref{bio1}  and $C(T_0)$
is a constant depending only of $T_0$.
\end{teor}

\demo Since it is similar to that of
Theorem 3.4 from \cite{ML}, we only give the main ideas. Let
$(\theta_m)_{m\in\mathbb{Z}^*}\subset
L^2\left(-\frac{\widetilde{T}}{2},\frac{\widetilde{T}}{2}\right)$
be the biorthogonal sequence from Theorem \ref{bio1}. For any $a>0$ define
$k_a=\frac{\sqrt{2\pi}}{a^2}(\chi_a*\chi_a)$, where $\chi_a$
represents the characteristic function $\chi_{[-a/2,a/2]}$.
Evidently, supp$(k_a)\subset [-a,a]$. We introduce
the function $\rho_m(x)=e^{ix\Im(\lambda_m)}k_a(x)$ and we define
\begin{equation}\label{zetam1}
\zeta_m=\frac{1}{\sqrt{2\pi}
\widehat{\rho}_m(i\overline{\lambda}_m)}\,\theta_m*\rho_m
\qquad(m\in\mathbb{Z}^*),
\end{equation}
where $\widehat{\rho}_m$ is the Fourier transform of $\rho_m$.
Evidently, $\zeta_m\in L^2\left(-\frac{\widetilde{T}}{2}-a,\frac{\widetilde{T}}{2}+a\right)$.
Let $T_0=\widetilde{T}+2a$. From the convolution's properties,
it follows that $(\zeta_m)_{m\in\mathbb{Z}^\ast}$ is a biorthogonal sequence
to the family $\left(e^{\lambda_nt}\right)_{n\in\mathbb{Z}^\ast}$ in
$L^2\left(-\frac{T_0}{2},\frac{T_0}{2}\right)$ and (\ref{est2}) is proved.\ctd

\section{Controllability results}

Now we are able to prove the main result of this paper.

{\em Proof of Theorem \ref{theo-main}:}
 Let $T>$max$\{2\beta,T_0\}$ and $(\zeta_m)_{m\in\mathbb{Z}^\ast}$ as in
Theorem \ref{theo-2}. We construct a control $v_\varepsilon\in
L^2(0,T)$ of \eqref{ec.in} corresponding to the initial data
$(u^0,u^1)\in {\cal H}_0$ as follows
\begin{equation}\label{eq.cont1}
v_\varepsilon(t)=-\sum_{m\in\mathbb{Z}^*}
 \frac{e^{-\overline{\lambda}_{m}\frac{T}{2}}}{\widehat{f}_{|m|}}
 \left(\widehat{u}^1_{|m|}+\lambda_{m}\widehat{u}^0_{|m|}\right)\widetilde{\zeta}_m\left( t-\frac{T}{2}
\right)\qquad (t\in(0,T)),
\end{equation}
where $\widetilde{\zeta}_m$ is the extension by zero of $\zeta_m$
to the interval $\left(-\frac{T}{2},\frac{T}{2}\right).$ From the
properties of the biorthogonal sequence
$(\zeta_{m})_{m\in\mathbb{Z}^\ast}$, it is easy to see that
$v_\varepsilon$ verifies (\ref{ec2.2}). To conclude that
$v_\varepsilon$ is a control for \eqref{ec.in}, we only have to
prove that the right hand side of (\ref{eq.cont1}) converges in
$L^2(0,T)$. This follows immediately from Theorem \ref{theo-2} and
the fact that $(u^0,u^1)\in {\cal H}_0$. Indeed, we have that
$$
\begin{array}{c}
\ds\int_0^T|v_{\ve}(t)|^2dt=\int_0^T\left|-
\sum_{m\in\mathbb{Z}^*}
 \frac{e^{-\overline{\lambda}_{m}\frac{T}{2}}}{\widehat{f}_{|m|}}
 \left(\widehat{u}^1_{|m|}+\lambda_{m}\widehat{u}^0_{|m|}\right)\widetilde{\zeta}_m\left( t-\frac{T}{2}
\right) \right|^2 dt=\\
\\
\ds=\int_{- \frac{\widetilde{T}}{2}-a}^{
\frac{\widetilde{T}}{2}+a}\left|-
\sum_{m\in\mathbb{Z}^*}
 \frac{e^{-\overline{\lambda}_{m}\frac{T}{2}}}{\widehat{f}_{|m|}}
 \left(\widehat{u}^1_{|m|}+\lambda_{m}\widehat{u}^0_{|m|}\right)\zeta_m\left(t\right)\right|^2 dt\leq C(
T_0)\|(u^0,u^1)\|^2_{\mathcal{H}_0}.
\end{array}
$$
The last inequality results from \eqref{est2} with the constant
$C(T_0)$ independent of $\varepsilon$ and $m$. Thus, the family of
controls $(v_\varepsilon)_{\varepsilon\in(0,1)}$ is uniformly
bounded in $L^2(0,T)$. In order to show that any weak limit of the
family $(v_\ve)_{\ve\in(0,1)}$ is a control for (\ref{ec.in0}) we
only have to pass to the limit as $\varepsilon$ goes to zero in
(\ref{ec2.2}).\ctd

\begin{obs}
Theorem \ref{theo-main} gives a uniform controllability result in
a time $T$ sufficiently large, but independent of $\ve$. From
Propositions \ref{lm1}, \ref{lmprod}  and \ref{lm20} we can give
an explicit expression of $T$. However, much more precise
estimates are needed in order to obtain the value $T=2\pi$ which
is known to be optimal in the limit case $\ve=0$.
\end{obs}

\section{Appendix}

The aim of this section is to give the proof of Lemma \ref{majqm}
from Subsection \ref{sec31}. Through this section $C$ denotes an
absolute positive constant.

{\em Proof of Lemma \ref{majqm}:} From the symmetry of the
sequence $(\lm_n)_{n\in\mathbb{Z}^*}$, it is sufficient to
consider only the case $m\in\mathbb{N}^{*}$. We have that

$
\begin{array}{c}
\ds\prod\limits_{\substack{n\in\mathbb{Z}^{*}\\n\neq m}}
\left|\frac{\lambda_{n}}{\lambda_{n}-\lambda_{m}}\right|^{2}=\ds\frac{m^2+\ve^2
m^{4\alpha}}{4m^2}\prod\limits_{\substack{n=1\\n\neq m}}^{\infty}
\frac{\ds (n^{2}+\ve^{2}n^{4\al})^{2}}
{\ds\left[|n-m|^{2}+\ve^{2}|n^{2\al}-m^{2\al}|^{2}
\right] \left[(n+m)^{2}+\ve^{2}|n^{2\al}-m^{2\al}|^{2}\right]}=\\
\\
=\ds\frac{m^2+\ve^2 m^{4\alpha}}{4m^2}
\underbrace{\prod\limits_{\substack{n=1\\n\neq m}}^{\infty}
\frac{|\lm_{n}|^{4}}{|\lm_{|n-m|}|^{2}|\lm_{m+n}|^{2}}}_{Q_{m}^1}
\underbrace{\prod\limits_{\substack{n=1\\n\neq m}}^{\infty}
\frac{\ds\left[|n-m|^{2}+\ve^{2}|n-m|^{4\al}\right]
\left[(n+m)^{2}+\ve^{2}(n+m)^{4\al}\right]} {\ds\left[|n-m|^{2}+
\ve^{2}|n^{2\al}-m^{2\al}|^{2}\right]
\left[(n+m)^{2}+\ve^{2}|n^{2\al}-m^{2\al}|^{2}\right]}}_{Q_{m}^2}.
\end{array}
$

Since, $$\ds Q_{m}^1=\prod\limits_{\substack{n=1\\n\neq m}}^
{\infty}\frac{|\lm_{n}|^{4}}{|\lm_{|n-m|}|^{2}|\lm_{m+n}|^{2}}\leq
\frac{\ds\prod_{\substack{n=1\\n\neq
m}}^{\infty}|\lambda_{n}|^{4}}
{\ds\prod_{n=1}^{m-1}|\lm_{n}|^{2}\prod_{n=1}^{\infty}|\lm_{n}|^{2}
\prod_{\substack{n=m+1\\n\neq2m}}^{\infty}|\lm_{n}|^{2}}\leq
\frac{|\lm_{2m}|^{2}}{|\lm_{m}|^{2}},$$ it follows that
\begin{equation}\label{eq0}
Q_{m}^1\leq16\qquad(m\in\mathbb{N}^{*}).
\end{equation}

For the product $Q_m^2$ we have that
$$
\begin{array}{l}
Q_{m}^2=\ds\prod\limits_{\substack{n=1\\n\neq m}}^{\infty}
\left(1+\frac
{|n-m|^{2}(n+m)^{2}\left[\ve^{2}m^{4\al-2}f(\frac{n}{m})+
\ve^{4}m^{8\al-4}g(\frac{n}{m})\right]}
{|n-m|^{2}(n+m)^{2}\left(1+
\ds\ve^{2}\frac{|n^{2\al}-m^{2\al}|^{2}}{|n-m|^{2}}\right)
\left(1+\ds\ve^{2}\frac{|n^{2\al}-m^{2\al}|^{2}}{(n+m)^{2}}\right)}\right)=\\
\\
=\exp\left(\ds\sum\limits_{\substack{n=1\\n\neq m}}^{\infty}
\ln\left(1+\frac {\ve^{2}m^{4\al-2}
f(\frac{n}{m})+\ve^{4}m^{8\al-4}g(\frac{n}{m})}{
\left(1+\ds\ve^{2}\frac{|n^{2\al}-m^{2\al}|^{2}}{|n-m|^{2}}\right)
\left(1+\ds\ve^{2}\frac{|n^{2\al}-m^{2\al}|^{2}}{(n+m)^{2}}\right)}\right)\right)\leq\\
\\
\leq\ds\exp\left(\sum\limits_{\substack{n=1\\n\neq m}}^{\infty}
\frac
{\ve^{2}m^{4\al-2}f(\frac{n}{m})+\ve^{4}m^{8\al-4}g(\frac{n}{m})}
{\left(1+
\ds\ve^{2}\frac{|n^{2\al}-m^{2\al}|^{2}}{|n-m|^{2}}\right)
\left(1+\ds\ve^{2}\frac{|n^{2\al}-m^{2\al}|^{2}}{(n+m)^{2}}
\right)}\right),
\end{array}
$$
where $f,g:(0,\infty)\to\mathbb{R}$ are the functions defined by
$$f(t)=|t-1|^{4\al-2}+(t+1)^{4\al-2} -2(t^{2}+1)\frac{(t^{2\al}-1)^{2}} {(t^{2}-1)^{2}},$$
$$g(t)=|t^{2}-1|^{4\al-2}-\frac{(t^{2\al}-1)^{4}} {(t^{2}-1)^{2}}.$$
Let us remark that, for any $\al\in[0,1)\setminus\left\{\frac{1}{2}\right\}$, the function $g$ can be bounded as follows
\begin{equation}\label{eq1}
g(t)\leq C\left\{
\begin{array}{ll}
t^{2\al}&t<\frac{1}{2}\\
|t^{2}-1|^{4\al-2}&t\in\left[\frac{1}{2},2\right]\\
t^{6\al-4}&t>2.
\end{array}\right.
\end{equation}

We prove the following inequality
\begin{equation}\label{eq2}
S_{g}:=\sum\limits_{\substack{n=1\\n\neq m}}^{\infty} \frac
{\ve^{4}m^{8\al-4}g(\frac{n}{m})} {\left(1+
\ds\ve^{2}\frac{|n^{2\al}-m^{2\al}|^{2}}{|n-m|^{2}}\right)
\left(1+\ds\ve^{2}\frac{|n^{2\al}-m^{2\al}|^{2}}{(n+m)^{2}}
\right)} \leq C |\Re(\lm_m)|\qquad(m\in\mathbb{N}^{*}).
\end{equation}
For that we write $S_g$ as follows
$$
\begin{array}{l}
S_{g}=\ds \left(\sum_{n=1}^{\left[\frac{m}{2}\right]}+
\sum\limits_{\substack{n=\left[\frac{m}{2}\right]+1\\n\neq
m}}^{2m-1}+\sum_{n=2m}^{\infty}\right) \frac {\ve^{4}m^{8\al-4}
g(\frac{n}{m})}
{\left(1+\ds\ve^{2}\frac{|n^{2\al}-m^{2\al}|^{2}}{|n-m|^{2}}\right)
\left(1+\ds\ve^{2}\frac{|n^{2\al}-m^{2\al}|^{2}}{(n+m)^{2}}\right)}
=S_{g}^{1}+S_{g}^{2}+S_{g}^{3}.
\end{array}
$$
In order to evaluate $S_g^i$ we use the following inequalities
\begin{equation}\label{inec1}
\frac{|n^{2\al}-m^{2\al}|^{2}}{|n-m|^{2}}\geq|n-m|^{4\al-2},\quad
\frac{|n^{2\al}-m^{2\al}|^{2}}{|n-m|^{2}}\geq(n+m)^{4\al-2},
\end{equation}
which hold for every $\al\in\left(\frac{1}{2},1\right)$  and
$n\in\mathbb{N}^*,\;n\neq m.$

For any
$\al\in\left[0,1\right)\setminus\left\{\frac{1}{2}\right\}$ we
evaluate the sum $S_{g}^1$ by using  \eqref{eq1} and by taking
into account that $m+n\leq2m$ and $m-n\geq\frac{1}{2}m$ for every
$n\in\left[1,\frac{m}{2}\right]$. We deduce that
$$
\ds S_{g}^{1}\leq C\sum_{n=1}^{\left[\frac{m}{2}\right]}\frac
{\ve^{4}m^{6\al-4}n^{2\al}} {\left(1+\ds\ve^{2}m^{4\al-2}\right)^2
}\leq \frac{C\ve^{4}m^{6\al-4}}{(1+\ve^{2}m^{4\al-2})^{2}}
\int_{1}^{\left[\frac{m}{2}\right]+1}t^{2\al}dt
\leq\frac{C\ve^{4}m^{8\al-3}}{(1+\ve^{2}m^{4\al-2})^{2}}\leq
C|\Re(\lm_m)|.
$$

Similarly, for any
$\al\in\left[0,1\right)\setminus\left\{\frac{1}{2}\right\}$ we
evaluate the sum $S_g^3$ as follows
$$
\begin{array}{c}
S_{g}^{3}\ds\leq\sum_{n=2m}^{\infty}\frac{C\ve^{4}m^{2\al}n^{6\al-4}}
{\left(1+\ds\ve^{2}\frac{|n^{2\al}-m^{2\al}|^{2}}{|n-m|^{2}}\right)
\left(1+\ds\ve^{2}\frac{|n^{2\al}-m^{2\al}|^{2}}{(n+m)^{2}}\right)}
\leq\ds\sum_{n=2m}^{\infty}\frac{C\ve^{4}m^{2\al}n^{6\al-4}}{(1+\ve^{2}n^{4\al-2})^{2}}
\leq C|\Re(\lm_m)|.
\end{array}
$$
For  $\alpha\in\left[0,\frac{1}{2}\right)$ we analyze $ S_{g}^{2}$ as follows
$$S_{g}^{2}\leq\ds\sum\limits_{n=\left[\frac{m}{2}\right]+1}^{2m-1}
\frac{C\ve^{4}|n^2-m^2|^{4\al-2}}{\left(1+\ds\ve^{2}
\frac{|n^{2\al}-m^{2\al}|^{2}}{|n-m|^{2}}\right)
\left(1+\ds\ve^{2}\frac{|n^{2\al}-m^{2\al}|^{2}}{(n+m)^{2}}\right)}\leq
\sum\limits_{n=\left[\frac{m}{2}\right]+1}^{2m-1}
C\ve^{4}m^{4\al-2}|n-m|^{4\al-2}\leq C|\Re(\lm_m)|.$$

And, for  $\al\in\left(\frac{1}{2},1\right)$, we have to
treat separately the cases $m\leq\gv$ and $m>\gv.$ For
$m\leq\gv$, we notice that the function $g$ is continuous  on a
compact set, so there exists a positive constant $C$ independent
of $\ve$ and $m$ such that $g(t)\leq C$. By using again
\eqref{inec1} it follows that
$$
\begin{array}{c}
S_{g}^{2}\leq\ds\sum\limits_{n=\left[\frac{m}{2}\right]+1}^{2m-1}
\frac{C\ve^{4}m^{8\al-4}}{\left(1+\ve^{2}|n-m|^{4\al-2}\right)
\left(1+\ds\ve^{2}\frac{|n-m|^{4\al}}{m^{2}}\right)}\leq
\sum_{k=1}^{m}\frac{C\ve^{4}m^{8\al-2}}{\left(1+\ve^{2}k^{4\al-2}\right)
\left(m^{2}+\ds\ve^{2}k^{4\al}\right)}\leq C|\Re(\lm_m)|.
\end{array}
$$
For the case $m>\gv$ it follows that
$$
\begin{array}{c}
S_{g}^{2}\leq\ds\sum_{n=\left[\frac{m}{2}\right]+1}^{2m-1}
\frac{C\ve^{4}m^{4\al-2}|n-m|^{4\al-2}} {\ve^{2}m^{4\al-2}
\left(1+\ds\ve^{2}\frac{(m+n)^{4\al-2}|n-m|^{2}}{m^{2}}\right)}\leq
\ds\sum_{k=1}^{m}\frac{C\ve^{2}m^{2}k^{4\al-2}} {m^{2}+\ds\ve^{2}m^{4\al-2}k^{2}}\leq\\
\\
\ds\leq
C\ve^{2}\sum_{k=1}^{\left[\frac{m^{2-2\al}}{\ve}\right]}k^{4\al-2}+
Cm^{4-4\al}
\sum_{k=\left[\frac{m^{2-2\al}}{\ve}\right]+1}^{m}k^{4\al-4}\leq
C|\Re(\lm_m)|,
\end{array}
$$
 which concludes the proof
of (\ref{eq2}).

Let us remark that, for any
$\alpha\in[0,1)\setminus\left\{\frac{1}{2}\right\}$, the function
$f$ can be bounded in the following way
\begin{equation}\label{eq3}
f(t)\leq C\left\{
\begin{array}{ll}
t^{2\al}&t<\frac{1}{2}\\
|t-1|^{4\al-2}&t\in\left[\frac{1}{2},2\right]\\
 t^{2\al-2}&t>2.
\end{array}
\right.
\end{equation}

We prove the following inequality
\begin{equation}\label{eq4}
S_{f}:= \sum\limits_{\substack{n=1\\n\neq m}}^{\infty} \frac
{\ve^{2}m^{4\al-2}f(\frac{n}{m})} {\left(1+
\ds\ve^{2}\frac{|n^{2\al}-m^{2\al}|^{2}}{|n-m|^{2}}\right)
\left(1+\ds\ve^{2}\frac{|n^{2\al}-m^{2\al}|^{2}}{(n+m)^{2}}
\right)}\leq C |\Re(\lm_m)|\qquad(m\in\mathbb{N}^{*}).
\end{equation}
Indeed, we have that
$$
\begin{array}{c}
S_{f}=\ds\left(\sum_{n=1}^{\left[\frac{m}{2}\right]}+
\sum\limits_{\substack{n=\left[\frac{m}{2}\right]+1\\n\neq
m}}^{2m-1}+\sum_{n=2m}^{\infty}\right) \frac
{C\ve^{2}m^{4\al-2}f(\frac{n}{m})}
{\left(1+\ds\ve^{2}\frac{|n^{2\al}-m^{2\al}|^{2}}{|n-m|^{2}}\right)
\left(1+\ds\ve^{2}\frac{|n^{2\al}-m^{2\al}|^{2}}{(n+m)^{2}}\right)}=S_{f}^{1}+S_{f}^{2}+S_f^3.
\end{array}
$$

For any
$\al\in\left[0,1\right)\setminus\left\{\frac{1}{2}\right\}$ we
evaluate the sum $S_{f}^1$ by taking into account \eqref{eq3} and
the fact that $m+n\leq2m$ and $m-n\geq\frac{1}{2}m$ for every
$n\in\left[1,\frac{m}{2}\right]$. We deduce that

$
S_{f}^{1}=\ds\sum\limits_{\substack{n=1}}^{\left[\frac{m}{2}\right]}\frac
{C\ve^{2}m^{2\al-2}n^{2\al}}
{\left(1+\ve^{2}\frac{(m^{2\al}-n^{2\al})^{2}}{(m-n)^{2}}\right)
\left(1+\ve^{2}\frac{(m^{2\al}-n^{2\al})^{2}}{(n+m)^{2}}\right)}\leq
\ds\frac{C\ve^{2}m^{2\al-2}}{\left(1+\ve^{2}m^{4\al-2}\right)^2}
\sum_{n=1}^{\left[\frac{m}{2}\right]}n^{2\al}\leq
 C|\Re(\lm_m)|.
$

Similarly, for any
$\al\in[0,1)\setminus\left\{\frac{1}{2}\right\}$, we deduce that
$S_f^3$ is bounded by $C|\Re(\lm_m)|$. Indeed,

$ S_{f}^3=\ds \sum\limits_{\substack{n=2m}}^{\infty} \frac
{\ve^{2}m^{2\al}n^{2\al-2}} {\left(1+
\ds\ve^{2}\frac{|n^{2\al}-m^{2\al}|^{2}}{|n-m|^{2}}\right)
\left(1+\ds\ve^{2}\frac{|n^{2\al}-m^{2\al}|^{2}}{(n+m)^{2}}
\right)}\leq \sum\limits_{\substack{n=2m}}^{\infty} \frac
{\ve^{2}m^{2\al}n^{2\al-2}}{\left(1+\ve^2n^{4\al-2}\right)^2}\leq
C|\Re(\lm_m)|.$

For $\al\in\left[0,\frac{1}{2}\right)$ we evaluate $S_f^2$, as follows
$$
\begin{array}{c}
\ds S_{f}^{2}=\sum_{n=\left[\frac{m}{2}\right]+1}^{2m-1} \frac
{C\ve^{2}|n-m|^{4\al-2}}
{\left(1+\ds\ve^{2}\frac{|n^{2\al}-m^{2\al}|^{2}}{|n-m|^{2}}\right)
\left(1+\ds\ve^{2}
\frac{|n^{2\al}-m^{2\al}|^{2}}{(n+m)^{2}}\right)}\leq
\sum_{n=\left[\frac{m}{2}\right]+1}^{2m-1} \frac
{C\ve^{2}|n-m|^{4\al-2}}
{1+\ds\ve^{2}\frac{|n^{2\al}-m^{2\al}|^{2}}{|n-m|^{2}}}\leq\\
\\
\ds\leq\sum_{n=\left[\frac{m}{2}\right]+1}^{2m-1} \frac
{C\ve^{2}|n-m|^{4\al}}
{|n-m|^{2}+\ds\ve^{2}m^{4\al}}=\sum_{k=1}^{m} \frac
{C\ve^{2}k^{4\al}}
{k^{2}+\ds\ve^{2}m^{4\al}}\leq\sum_{k=1}^{m} \frac
{C\ve^{2}m^{4\al}}
{k^{2}+\ds\ve^{2}m^{4\al}}\leq C|\Re(\lm_m)|.
\end{array}
$$
If $\al\in\left(\frac{1}{2},1\right)$ and $m\leq\gv$, we have that
$$
S_{f}^{2}\leq\ds\sum_{n=\left[\frac{m}{2}\right]+1}^{2m-1} \frac {C\ve^{2}|n-m|^{4\al-2}m^2}
{\left(1+\ve^{2}m^{4\al-2}\right)\left(m^{2}+\ve^{2}|n-m|^{2}m^{4\al-2}\right)}\leq
C\ve^{2}m^{2}\sum_{k=1}^{m}k^{4\al-2}\leq C|\Re(\lm_m)|,
$$
and for $\al\in\left(\frac{1}{2},1\right)$ and $m>\gv$ the following estimates takes place
$$
\begin{array}{c}
\ds S_{f}^{2}\leq\ds\sum_{n=\left[\frac{m}{2}\right]+1}^{2m-1} \frac {C\ve^{2}|n-m|^{4\al-2}m^2}
{\left(1+\ve^{2}m^{4\al-2}\right)\left(m^{2}+\ve^{2}|n-m|^{2}m^{4\al-2}\right)}\leq
\frac{C\ve^{2}m^{2}}{1+\ve^{2}m^{4\al-2}}
\sum_{k=1}^{m}\frac{k^{4\al-2}}{m^{2}+\ve^{2}k^{2}m^{4\al-2}}\leq\\
\\
\ds \leq C m^{4-4\al}\left(\sum_{k=1}^{\left[\frac{m^{2-2\al}}{\ve}\right]}\frac{k^{4\al-2}}{m^2}+
\sum_{k=\left[\frac{m^{2-2\al}}{\ve}\right]+1}^m\frac{k^{4\al-4}}{\ve^2m^{4\al-2}}\right)\leq C|\Re(\lm_m)|.
\end{array}
$$
Now from (\ref{eq0}), (\ref{eq2}) and (\ref{eq4}) it results
(\ref{eq00}) and the proof of Lemma \ref{majqm} ends.\ctd

\
\

\textbf{Acknowledgement} The first author was partially supported
by the strategic grant POSDRU/CPP107/ DMI1.5/S/78421, Project ID
78421 (2010), co-financed by the European Social Fund - Investing
in People, within the Sectorial Operational Programme Human
Resources Development 2007-2013. The second author was partially
supported by Grant PN-II-ID-PCE-2011-3-0257 of the Romanian
National Authority for Scientific Research, CNCS – UEFISCDI  and
by Grant MTM2011-29306 funded by MICINN (Spain).

\end{document}